\documentclass[11pt,reqno]{amsart}
\usepackage{cmap}
\usepackage[T1]{fontenc}
\usepackage[utf8]{inputenc}
\usepackage[english]{babel}
\usepackage{amsfonts,amssymb,amsthm,amsmath}
\usepackage{url}
\usepackage{enumitem}

\usepackage{secdot}

\allowdisplaybreaks

\usepackage{graphicx}
\usepackage{epstopdf}
\usepackage{subfig}

\usepackage{a4wide}
\usepackage{xcolor}
\usepackage[colorlinks=true,linktocpage,pdfpagelabels,
bookmarksnumbered,bookmarksopen]{hyperref}
\definecolor{ForestGreen}{rgb}{0.1,0.6,0.05}
\definecolor{EgyptBlue}{rgb}{0.063,0.1,0.6}
\hypersetup{
	colorlinks=true,
	linkcolor=EgyptBlue,         
	citecolor=ForestGreen,
	urlcolor=olive
}

\usepackage[hyperpageref]{backref}

\newtheorem{theorem}{Theorem}

\newtheorem{corollary}[theorem]{Corollary}

\theoremstyle{definition}
\newtheorem{definition}[theorem]{Definition}
\newtheorem{remark}[theorem]{Remark}

\let\OLDthebibliography\thebibliography
\renewcommand\thebibliography[1]{
	\OLDthebibliography{#1}
	\setlength{\parskip}{1pt}
	\setlength{\itemsep}{1pt plus 0.3ex}
}

\numberwithin{equation}{section}
\numberwithin{theorem}{section}
\numberwithin{equation}{section}
\numberwithin{theorem}{section}

\newcommand{\N}{\mathbb{N}}
\newcommand{\R}{\mathbb{R}}

\title[Basis properties of Fu\v{c}\'ik eigenfunctions]{Basis properties of Fu\v{c}\'ik eigenfunctions}

\author[F.~Baustian]{Falko Baustian}
\author[V.~Bobkov]{Vladimir Bobkov}

\address[F.~Baustian]{\newline\indent
	Department of Mathematics, University of Rostock,
	\newline\indent 
	Ulmenstra{\ss}e 69, 18057 Rostock, Germany
}
\email{falko.baustian@uni-rostock.de}

\address[V.~Bobkov]{\newline\indent
	Department of Mathematics and NTIS, Faculty of Applied Sciences,
	\newline\indent 
	University of West Bohemia, Univerzitn\'i 8, 301 00 Plze\v{n}, Czech Republic\\
	\newline\indent 
	Institute of Mathematics, Ufa Federal Research Centre, RAS,
	\newline\indent 
	Chernyshevsky str. 112, 450008 Ufa, Russia
}
\email{bobkov@kma.zcu.cz}

\date{}


\thanks{
	V.~Bobkov was supported by the project LO1506 of the Czech Ministry of Education, Youth and Sports, and by the grant 18-03253S of the Grant Agency of the Czech Republic.
}

\begin{document}
\begin{abstract} 
	We establish sufficient assumptions on sequences of Fu\v{c}\'ik eigenvalues of the one-dimensional Laplacian which guarantee that the corresponding Fu\v{c}\'ik eigenfunctions form a Riesz basis in $L^2(0,\pi)$.
\end{abstract}
\maketitle 
	
\section{Introduction}

The classical spectral theorem for compact self-adjoint (linear) operators on a Hilbert space $X$ asserts the existence of a sequence of eigenfunctions of such an operator that forms an orthogonal basis of $X$, see, e.g., \cite{RS}. 
One of the simplest examples is given by the sine functions $\sin(nx)$, $n \in \mathbb{N}$, which are eigenfunctions of the one-dimensional Dirichlet Laplacian, i.e.,
$$
-u'' = \lambda u 
~~ \text{in}~ (0,\pi),
\quad u(0)=u(\pi)=0,
$$
and hence they form an orthogonal basis in $L^2(0,\pi)$. 
Although the assumptions of the spectral theorem can not be weakened in general, eigenfunctions of certain operators which do not satisfy the imposed requirements might still form a basis of a space $X$. 
Classes of such operators have been an active topic of research, and a great amount of significant results have been obtained. 
We refer to \cite{agran,djak,kats,locker,markus,skal1} and the extensive bibliographies therein for a deeper discussion.

The references indicated above are concerned mainly with linear operators. 
At the same time, the conclusion of the spectral theorem can be valid even for some nonlinear operators.
As an example, there exists $p_0>1$ such that the generalized trigonometric functions $\sin_p(nx)$ defined as eigenfunctions
of the one-dimensional Dirichlet $p$-Laplacian with $p > p_0$, i.e.,
$$
-(|u'|^{p-2}u')' = \lambda |u|^{p-2} u 
~~ \text{in}~ (0,\pi_p),
\quad u(0)=u(\pi_p)=0,
$$
form a Riesz basis in $L^2(0,\pi_p)$, see \cite{bind}, and \cite{BM,BE,EGL} for further developments.
Here, $\pi_p=\frac{2\pi}{p \sin(\pi/p)}$.
While this basis is non-orthogonal except for $p=2$, it has certain advantages, in particular, in a numerical study of nonlinear equations, see \cite{BL}.
Let us remark that 
it remains unknown whether the same basisness result holds true for all $p >1$.
Notice that the generalized trigonometric functions $\sin_p(nx)$ share with $\sin(nx)$ the anti-periodic structure which is significantly used in the proof of \cite{bind}. 
In general, there exist several results on assumptions for a function $f$ to guarantee that the corresponding system of dilated functions, i.e., the sequence of the form $\{f(nx)\}$, is a basis, see, e.g., \cite{BB,BM,HLS,HLS2,wintner2}, to mention a few.

In the present work, we investigate the basisness of eigenfunctions for another type of eigenvalue problems for the classical Dirichlet Laplacian in one dimension which does not fit into the framework of the  spectral theorem, namely,
\begin{equation}\label{eq:0}
	-u'' = \alpha u^{+}(x)-\beta u^{-}(x)
	~~ \text{in}~ (0,\pi),
	\quad u(0)=u(\pi)=0,
\end{equation}
where both $\alpha,\beta$ are spectral parameters.
This so-called \textit{Fu\v{c}\'ik} (or \textit{Dancer-Fu\v{c}\'ik}) \textit{eigenvalue problem} was originated in the works \cite{dancer} and \cite{fucik} in the context of studies of elliptic equations with ``jumping'' nonlinearities. 
Various aspects of the Fu\v{c}\'ik eigenvalue problem for second-order elliptic operators where intensively studied afterwards.
We refer to \cite{cac,schech2} for the local existence of curves in the Fu\v{c}\'ik spectrum emanated from classical eigenvalues $\lambda_k$, 
to \cite{castro,drabek} for variational characterizations of such curves, 
to \cite{pinsal,rynne} for their asymptotic behaviour, 
to \cite{cues,li} for nodal properties of the corresponding Fu\v{c}\'ik eigenfunctions,
and to \cite{dancer2,perera} for algebro-topological properties of the associated energy functional.
Although this list is far from being complete, many other relevant references can be found in the indicated papers.
Nevertheless, the structure and properties of the Fu\v{c}\'ik spectrum and Fu\v{c}\'ik eigenfunctions are understood only to a rather limited extent, and their investigation can be a hard task even in the case of operators acting in finite-dimensional spaces, see the discussion in \cite{holubova,loos}.
In particular, we are not aware of previous results on the basis properties of Fu\v{c}\'ik eigenfunctions of the problem \eqref{eq:0} for the one-dimensional Laplacian.

Fu\v{c}\'ik eigenfunctions associated with Fu\v{c}\'ik eigenvalues $(\alpha,\beta)$ of the problem \eqref{eq:0} can be interpreted as specific alterations of the sine functions that, in general, do not lead to a system of dilated or orthogonal functions when $\alpha \neq \beta$, although they preserve a certain periodic structure.
By establishing explicit formulas for the distances between Fu\v{c}\'ik eigenfunctions and corresponding sine functions, and making use of the periodic structure, we obtain sequences of Fu\v{c}\'ik eigenvalues $(\alpha(n),\beta(n))$ whose corresponding Fu\v{c}\'ik eigenfunctions form a Riesz basis in $L^2(0,\pi)$. 
To this end, we employ several classical results from the stability theory of Paley and Wiener \cite{paley}, which provides basic information on the basisness of  
Fu\v{c}\'ik eigenfunctions, and we leave more developed approaches for further investigation.

Let us outline the structure of this article. 
In the remainder of the introduction, we rigorously introduce normalized Fu\v{c}\'ik eigenfunctions and Fu\v{c}\'ik systems as sequences of such functions in Section \ref{sec:fucik}, and we state our main results in Section \ref{sec:main-results}. 
We derive explicit formulas for the norms of the normalized Fu\v{c}\'ik eigenfunctions, their distances to the corresponding sine functions, and some scalar products between Fu\v{c}\'ik eigenfunctions and sine functions in Section \ref{sec:norms}. 
In Section \ref{sec:asympotics}, we establish upper bounds on the distances which are essential for the proof of Theorem \ref{th:main}. 
The proof of Theorem \ref{thm:even} is given in Section \ref{sec:even}, and the proof of Theorem \ref{th:main2} is given in Section \ref{sec:thm2}. 
We conclude the article with final remarks in Section \ref{sec:final}.
Appendices \ref{sec:appendix} and \ref{sec:appendix2} contain several auxiliary technical details.

\subsection{\texorpdfstring{Fu\v{c}\'ik}{Fucik} eigenvalues and eigenfunctions}\label{sec:fucik}

\begin{figure}[h!]
	\center{\includegraphics[width=0.32\linewidth]{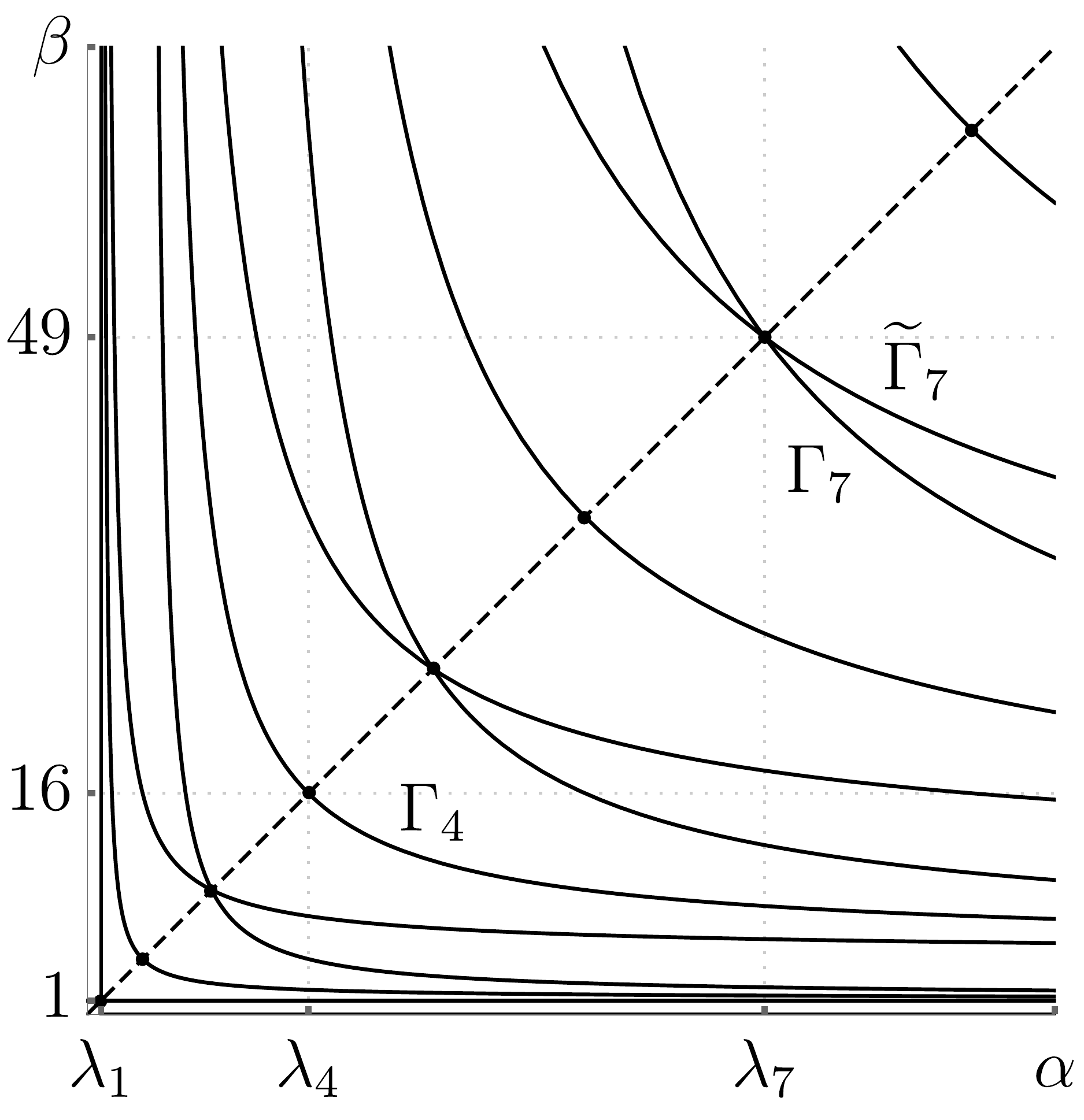}}
	\caption{Several curves of the Fu\v{c}\'ik spectrum}
	\label{fig1}
\end{figure}

The \textit{Fu\v{c}\'ik spectrum} $\Sigma(0,\pi) \subset \R^2$ 
of the linear Dirichlet Laplacian in one dimension is defined as the set of all pairs of parameters $(\alpha,\beta)\in\R^2$ for which the problem
\begin{equation}\label{eq:fucik}
\tag{\ref{eq:0}}
\left\{
\begin{alignedat}{1}
-u''(x) &= \alpha u^{+}(x)-\beta u^{-}(x) ~~ \mbox{in }(0,\pi), \\
u(0) &=u(\pi)=0,
\end{alignedat}
\right.
\end{equation}
has a non-trivial solution. 
Here, $u^{+}=\max(u,0)$ and $u^{-}=\max(-u,0)$ denote the positive and negative parts of $u$, respectively. 
Clearly, for $\alpha=\beta =: \lambda$ we obtain the standard eigenvalue problem $-u''= \lambda u$ with zero Dirichlet boundary conditions,
which possesses the complete sequence of eigenvalues $\lambda_n=n^2$, $n\in\N$, and the sine functions $\phi_n=\sin(nx)$ as corresponding eigenfunctions. Since $\phi_1$ has definite sign, the lines $\{1\}\times\R$ and $\R\times\{1\}$ are the trivial part of $\Sigma(0,\pi)$ and we have
$$
\Sigma(0,\pi)\setminus\left((\{1\}\times\R)\cup(\R\times\{1\})\right)
\subset
\{(\alpha,\beta)\in\R^2: \, \alpha>1, \,\beta>1\}
$$
by the variational characterization of the first eigenvalue $\lambda_1=1$.
Any $(\alpha,\beta) \in \Sigma(0,\pi)$ is called \textit{Fu\v{c}\'ik eigenvalue} and any corresponding non-zero solution of \eqref{eq:fucik} is called \textit{Fu\v{c}\'ik eigenfunction}.
Hereinafter, under a solution of \eqref{eq:fucik} we mean the classical (i.e., $C^2$-regular) solution of that problem.
Any Fu\v{c}\'ik eigenfunction consists of positive bumps of length $l_1=\frac{\pi}{\sqrt{\alpha}}$ alternating with negative bumps of length $l_2=\frac{\pi}{\sqrt{\beta}}$, and has the form of translations and multiples of $\sin(\sqrt{\alpha}x)$ at positive bumps and $\sin(\sqrt{\beta}x)$ at negative bumps.
For each even number $n$ of bumps there exists one curve
$$
\Gamma_n
=
\left\{(\alpha,\beta)\in\R^2: \,
\frac{n}{2}\frac{\pi}{\sqrt{\alpha}}
+
\frac{n}{2}\frac{\pi}{\sqrt{\beta}}
=\pi
\right\}
$$
in the Fu\v{c}\'ik spectrum $\Sigma(0,\pi)$ that contains $(\lambda_n,\lambda_n)$, while for each odd number $n\geq3$ of bumps there are two curves
\begin{align*}
	\Gamma_n &= \left\{(\alpha,\beta)\in\R^2: \,
	\frac{n+1}{2}\frac{\pi}{\sqrt{\alpha}}+\frac{n-1}{2}\frac{\pi}{\sqrt{\beta}}=\pi
	\right\},
	\\
	\widetilde{\Gamma}_n &=
	\left\{(\alpha,\beta)\in\R^2: \, \frac{n-1}{2}\frac{\pi}{\sqrt{\alpha}}+\frac{n+1}{2}\frac{\pi}{\sqrt{\beta}}=\pi
	\right\}
\end{align*}
in $\Sigma(0,\pi)$ intersecting at $(\lambda_n,\lambda_n)$, see Figure \ref{fig1}.
These curves completely describe the non-trivial part of $\Sigma(0,\pi)$, see \cite[Lemma 2.8]{fucik}.
Observe that if $u$ is a Fu\v{c}\'ik eigenfunction for the pair of parameters $(\alpha,\beta)$, then so is $tu$ for any $t>0$, while $-tu$ is a Fu\v{c}\'ik eigenfunction for $(\beta,\alpha)$. 
As a consequence, since $(\alpha,\beta) \in \Gamma_n$ for odd $n$ implies $(\beta,\alpha) \in \widetilde{\Gamma}_n$, we will neglect the curve $\widetilde{\Gamma}_n$ from our further consideration. 
Moreover, notice that if $u$ is a Fu\v{c}\'ik eigenfunction for the pair of parameters $(\alpha,\beta) \in \Gamma_{n}$ for even $n$, then $v(x)=u(\pi-x)$ is also a Fu\v{c}\'ik eigenfunction for $(\alpha,\beta)$ and $v \neq u$.

\begin{figure}[h!]
	\centering
	\subfloat[][\centering $n=2$]
	{\includegraphics[width=0.39\linewidth]{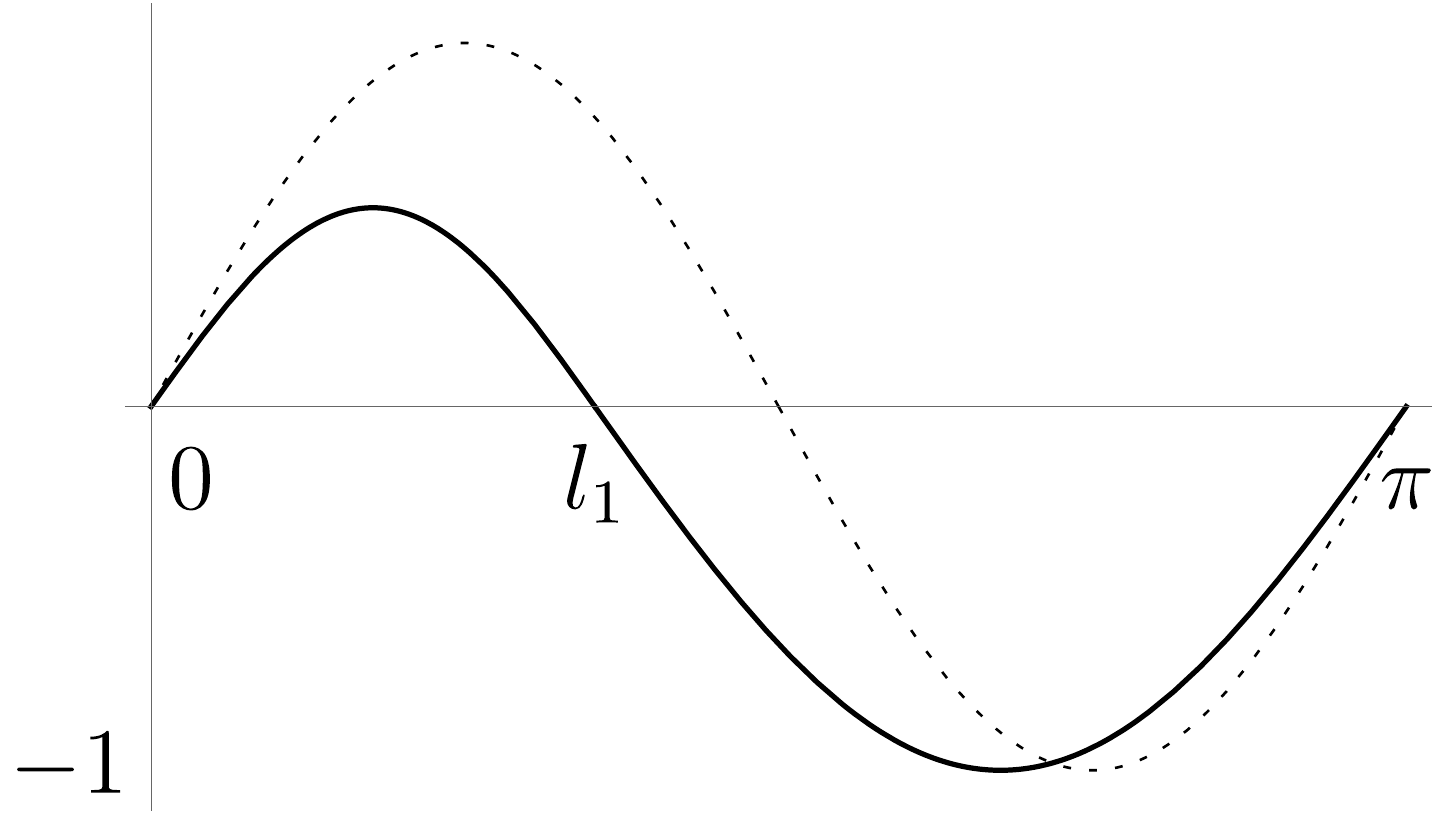}}
	\qquad\qquad
	\subfloat[][\centering $n=3$]
	{\includegraphics[width=0.39\linewidth]{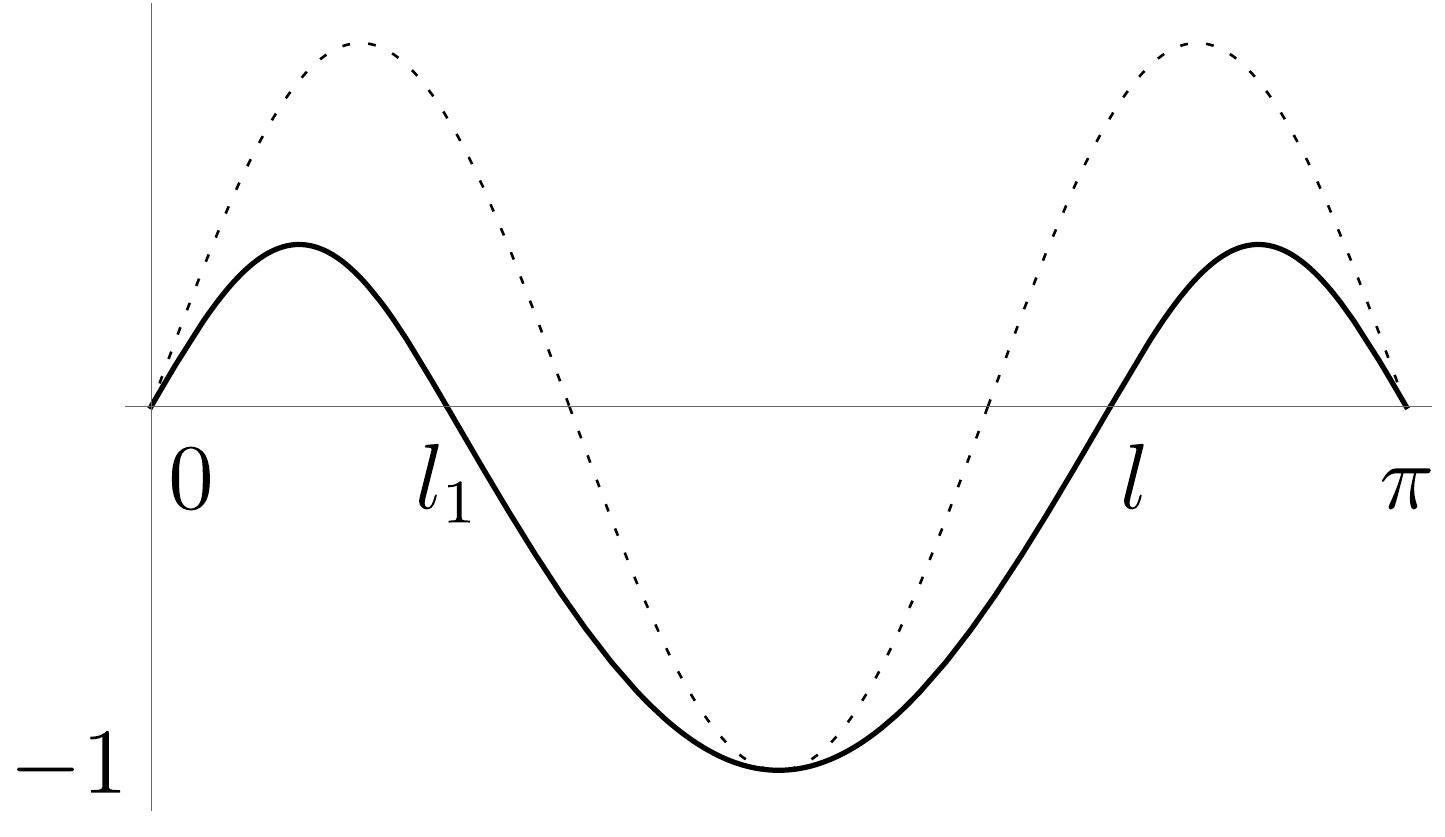}}
	\caption{A normalized Fu\v{c}\'ik eigenfunction $f^n_{\alpha,\beta}$ (solid) and $\sin(nx)$ (dotted)}
	\label{fig2}
\end{figure}

In order to uniquely specify a Fu\v{c}\'ik eigenfunction for each point of $\Sigma(0,\pi)$, let us introduce the following special choice of Fu\v{c}\'ik eigenfunctions, see Figure \ref{fig2}.
\begin{definition}
Let $n\geq2$ and $(\alpha,\beta)\in \Gamma_n$. The \textit{normalized Fu\v{c}\'ik eigenfunction} $f^n_{\alpha,\beta}$
is the $C^2$-solution of the boundary value problem \eqref{eq:fucik} with $(f^n_{\alpha,\beta})'(0)>0$  which is normalized by
$$
\|f^n_{\alpha,\beta}\|_{\infty}
=
\sup_{x\in[0,\pi]}|f^n_{\alpha,\beta}(x)|=1.
$$
For $n=1$, we set $f^1_{\alpha,\beta}=\phi_1$ for every $(\alpha,\beta)\in(\{1\}\times\R)\cup(\R\times\{1\})$.
\end{definition}

The normalized Fu\v{c}\'ik eigenfunctions can be described more explicitly by the following piecewise definition.
Let $n \geq 2$.
For $\alpha\geq n^2 \geq \beta$ we have
\begin{equation}\label{eq:fucikpiecewise1}
f^n_{\alpha,\beta}(x)=
\left\{
\begin{array}{clrl}
\frac{\sqrt{\beta}}{\sqrt{\alpha}}\sin(\sqrt{\alpha}\,(x-kl)) \quad &\mbox{for} 
&kl&\leq x<kl+l_1, \\
-\sin(\sqrt{\beta}\,(x-kl-l_1)) \quad &\mbox{for}
&kl+l_1&\leq x<(k+1)l,
\end{array}
\right.
\end{equation}
and for $\beta > n^2 > \alpha$ we have
\begin{equation}\label{eq:fucikpiecewise2}
f^n_{\alpha,\beta}(x)=
\left\{
\begin{array}{clrl}
\sin(\sqrt{\alpha}\,(x-kl)) \quad &\mbox{for} 
&kl&\leq x<kl+l_1, \\
-\frac{\sqrt{\alpha}}{\sqrt{\beta}}\sin(\sqrt{\beta}\,(x-kl-l_1)) \quad &\mbox{for}
&kl+l_1&\leq x<(k+1)l,
\end{array}
\right.
\end{equation}
where $l=l_1+l_2$ and $k \in \mathbb{N}_0 = \mathbb{N} \cup \{0\}$.
Notice that \eqref{eq:fucikpiecewise1} and \eqref{eq:fucikpiecewise2} define $f^n_{\alpha,\beta}$ on the whole $\mathbb{R}_+$. 
We also remark that $f^n_{\alpha,\beta} \not\in C^3[0,\pi]$ provided $\alpha \neq \beta$.

In our purpose to form a basis in $L^2(0,\pi)$ that consists solely of normalized Fu\v{c}\'ik eigenfunctions, we pick one normalized Fu\v{c}\'ik eigenfunction from each curve $\Gamma_n$, $n\geq2$, and the sine function as the normalized Fu\v{c}\'ik eigenfunction of the trivial part of $\Sigma(0,\pi)$.
In this regard, the following definition will be useful.
\begin{definition}\label{def:FS}
We define a \textit{Fu\v{c}\'ik system} $F_{\alpha,\beta}=\{f^n_{\alpha(n),\beta(n)}\}$ as a sequence of normalized Fu\v{c}\'ik eigenfunctions with the mappings $\alpha,\beta\colon\N\to\R$ satisfying $\alpha(1)=\beta(1)=1$ and $(\alpha(n),\beta(n))\in \Gamma_n$ for every $n\geq2$.
\end{definition}

\subsection{Main results}\label{sec:main-results}
In this section, we summarize our main results for the basisness of Fu\v{c}\'ik systems given by Definition \ref{def:FS}. 
We present growth conditions on the mappings $\alpha,\beta$ which guarantee that the Fu\v{c}\'ik system $F_{\alpha,\beta}$ is a Riesz basis in $L^2(0,\pi)$. 
We call a complete system $\{\psi_n\}$ in $L^2(0,\pi)$ a \textit{Riesz basis} if there exist positive constants $c,C>0$ such that the inequalities
$$
c\sum_{n=1}^N|\alpha_n|^2\leq\big\|\sum_{n=1}^N\alpha_n\psi_n\big\|\leq C\sum_{n=1}^N|\alpha_n|^2
$$
are satisfied for arbitrary $N\in\N$ and any constants $\alpha_1,\dots,\alpha_N$. 
In fact, Riesz bases are images of an orthonormal basis under a linear homeomorphism. Several equivalent characterizations of a Riesz basis can be found in \cite[Theorem 9]{young}.
Hereinafter, $\|\cdot\|$ denotes the standard norm in $L^2(0,\pi)$.

We make use of several methods from the stability theory of Paley and Wiener to show that a Fu\v{c}\'ik system $F_{\alpha,\beta}$ inherits the basis properties from the system of sine functions $\{\phi_n\}$ provided that the two sequences are sufficiently close to each other. 
There exist various concepts of nearness between systems of functions, among which we will be interested in the following two most classical notions, see, e.g., \cite{sing}.
\begin{definition}
	Let $\{\varphi_n\}$ and $\{\psi_n\}$ be two sequences of functions. The sequence $\{\psi_n\}$ is \textit{quadratically near} to $\{\varphi_n\}$ if 
	\begin{equation}\label{eq:quad_near}
	\sum_{n=1}^{\infty} \|\varphi_n-\psi_n\|^2=r<\infty
	\end{equation}
	for a constant $r \geq 0$, and \textit{strongly quadratically near} if \eqref{eq:quad_near} holds for $r<1$.
	The sequence $\{\psi_n\}$ is \textit{Paley-Wiener near} to $\{\varphi_n\}$ if there exits a constant $\lambda \in (0,1)$ such that
	\begin{equation*}\label{eq:palwie_near}
	\big\|\sum_{n=1}^{N} \alpha_n(\varphi_n-\psi_n)\big\|<\lambda\big\|\sum_{n=1}^{N}\alpha_n\varphi_n\big\|
	\end{equation*}
	holds for arbitrary $N\in\N$ and any constants $\alpha_1,\dots,\alpha_N$.
\end{definition}

We introduce our results, each one connected to one of these nearness concepts. The first general result, Theorem \ref{th:main}, allows a rather free choice of the Fu\v{c}\'ik system $F_{\alpha,\beta}$ and utilizes the strong quadratic nearness of that system to the sine functions $\{\phi_n\}$. 
In Theorem \ref{thm:even}, we significantly improve the constants of Theorem \ref{th:main} by means of the quadratic nearness of special Fu\v{c}\'ik systems $F_{\alpha,\beta}$  with $\alpha(n)=\beta(n)=n^2$ for odd $n$ to $\{\phi_n\}$.
For the final result the normalized Fu\v{c}\'ik eigenfunctions $f^n_{\alpha,\beta}$ for even $n$ are chosen in such a way that they form a sequence of dilated functions, while for odd $n$ we just pick the sine functions as in the previous case. 
This specific choice of the Fu\v{c}\'ik system allows us to apply the separation of variables approach of Duffin and Eachus \cite{duff} in order to establish the Paley-Wiener nearness to the sine functions. 
In view of the nature of Riesz bases, both approaches are intrinsically based on the construction of a bounded invertible operator $T: L^2(0,\pi) \to L^2(0,\pi)$ which maps the trigonometric system $\{\phi_n\}$ to the Fu\v{c}\'ik system $F_{\alpha,\beta}$.

The basisness of systems that are quadratically near to a complete orthonormal system was first studied by Bary in \cite{bary2}. 
A system $\{\psi_n\}$ which is quadratically near to a complete orthonormal system $\{\varphi_n\}$ is a Riesz basis provided that it is $\omega$-linearly independent, i.e., if the strong convergence
$$
\sum_{n=1}^{\infty} \eta_n\psi_n
= 
\lim_{m \to \infty} \|\sum_{n=1}^{m}\eta_n\psi_n\|=0
$$
for a sequence of scalars $\{\eta_n\}$ implies $\eta_n=0$ for every $n\in\N$. 
A proof of this stability result by means of compact operators is given, e.g., in \cite[Theorem V-2.20]{kato}. 
If the system $\{\psi_n\}$ satisfies the more restrictive assumption of being strongly quadratically near to a complete orthonormal system $\{\varphi_n\}$, then $\{\psi_n\}$  is also a Riesz basis, see, e.g., \cite[Corollary V-2.22]{kato}.
We establish a summation criterium for bounds on the mappings $\alpha$ and $\beta$ of the Fu\v{c}\'ik system $F_{\alpha,\beta}$ that yields the basisness of that system by means of the strong quadratic nearness to the system of sine functions $\{\phi_n\}$.

\begin{theorem}\label{th:main}
Let $F_{\alpha,\beta}$ be a Fu\v{c}\'ik system. 
For any natural $n \geq 2$, we set
\begin{equation}\label{eq:Bn}
C_n(x,y) =
\left\{
\begin{aligned}
&\frac{4(3+\pi^2)\pi}{9}\left(\frac{\max(\sqrt{x},\sqrt{y})}{n}-1\right)^2  &&\mbox{for even } 
n,&\\
&4\pi\frac{n^2(n^2+1)}{(n-1)^4}\left(\frac{\sqrt{x}}{n}-1\right)^2  &&\mbox{for odd } 
n\mbox{ with } x\geq y, \\
&5\pi\frac{n^2(n^2+1)}{(n+1)^4}\left(\frac{\sqrt{y}}{n}-1\right)^2  &&\mbox{for odd } 
n\mbox{ with } y>x.
\end{aligned}
\right.
\end{equation}
If the summation formula
\begin{equation}\label{eq:sumBn}
\sum_{n=2}^{\infty} C_n(\alpha(n),\beta(n))<\frac{\pi}{2}
\end{equation}
is satisfied, then $F_{\alpha,\beta}$ is a Riesz basis in $L^2(0,\pi)$.
\end{theorem}
The definition \eqref{eq:Bn} of $C_n$ is given by the bounds 
\eqref{eq:asymp:dist-fn-odd}, \eqref{eq:asymp-fn-odd-1}, \eqref{eq:asymp:dist-fn-odd2}, \eqref{eq:asymp:dist-fn-odd-b>a}
on the distances $\|f^n_{\alpha,\beta}-\phi_n\|^2$ that we will derive in Section \ref{sec:asympotics} below. 
The summation formula \eqref{eq:sumBn} guarantees that the Fu\v{c}\'ik system $F_{\alpha,\beta}$ is quadratically near to the system of sine functions $\{\phi_n\}$ 
in the sense that 
\begin{equation}\label{eq:sqna}
\sum_{n=1}^\infty
\|f^n_{\alpha,\beta}-\phi_n\|^2 < \frac{\pi}{2}.
\end{equation}
This inequality implies that the rescaled system $\sqrt{2/\pi} \, F_{\alpha,\beta}$ is strongly quadratically near to the complete orthonormal system  $\{\sqrt{2/\pi} \, \phi_n\}$ in $L^2(0,\pi)$, and thus it is a Riesz basis by \cite[Corollary V-2.22]{kato}. 
Hence, the initial Fu\v{c}\'ik system $F_{\alpha,\beta}$ is also a Riesz basis in $L^2(0,\pi)$.

We can make use of Theorem \ref{th:main} to give asymptotic bounds on the mappings $\alpha$ and $\beta$. 
\begin{corollary}\label{cor:thm1}
	Let $F_{\alpha,\beta}$ be a Fu\v{c}\'ik system and let $\varepsilon>0$ be fixed. 
	Let the mappings $\alpha$ and $\beta$ satisfy
	\begin{equation}\label{eq:maxalpha}
	\max\left(\sqrt{\alpha(n)},\sqrt{\beta(n)}\right)\leq n+\sqrt{c_n}\, n^{(1-\varepsilon)/2}
	\end{equation}
	for every $n\geq 2$ with non-negative constants
	\begin{align}
	\notag
	c_n&<\frac{9}{8(3+\pi^2)}\cdot\frac1{\zeta(1+\varepsilon)-1} \qquad~~ \mbox{for even }n, \\
	c_n&<\frac{(n-1)^4}{8n^2(n^2+1)}\cdot\frac1{\zeta(1+\varepsilon)-1} \quad~~ \mbox{for odd }n\mbox{ with }\alpha(n)\geq\beta(n), \label{eq:bas_con2} \\
	c_n&<\frac{(n+1)^4}{10n^2(n^2+1)}\cdot\frac1{\zeta(1+\varepsilon)-1} \quad \mbox{for odd }n\mbox{ with }\beta(n)>\alpha(n), \label{eq:bas_con3}
	\end{align}
	where $\zeta$ is the Riemann zeta function defined by $\zeta(s)=\sum_{n=1}^{\infty}\frac1{n^s}$. 
	Then $F_{\alpha,\beta}$ is a Riesz basis in $L^2(0,\pi)$.
\end{corollary}

\begin{figure}[h!]
	\centering
	\subfloat[\centering $\epsilon = 0.1$]
	{{\includegraphics[width=0.32\linewidth]{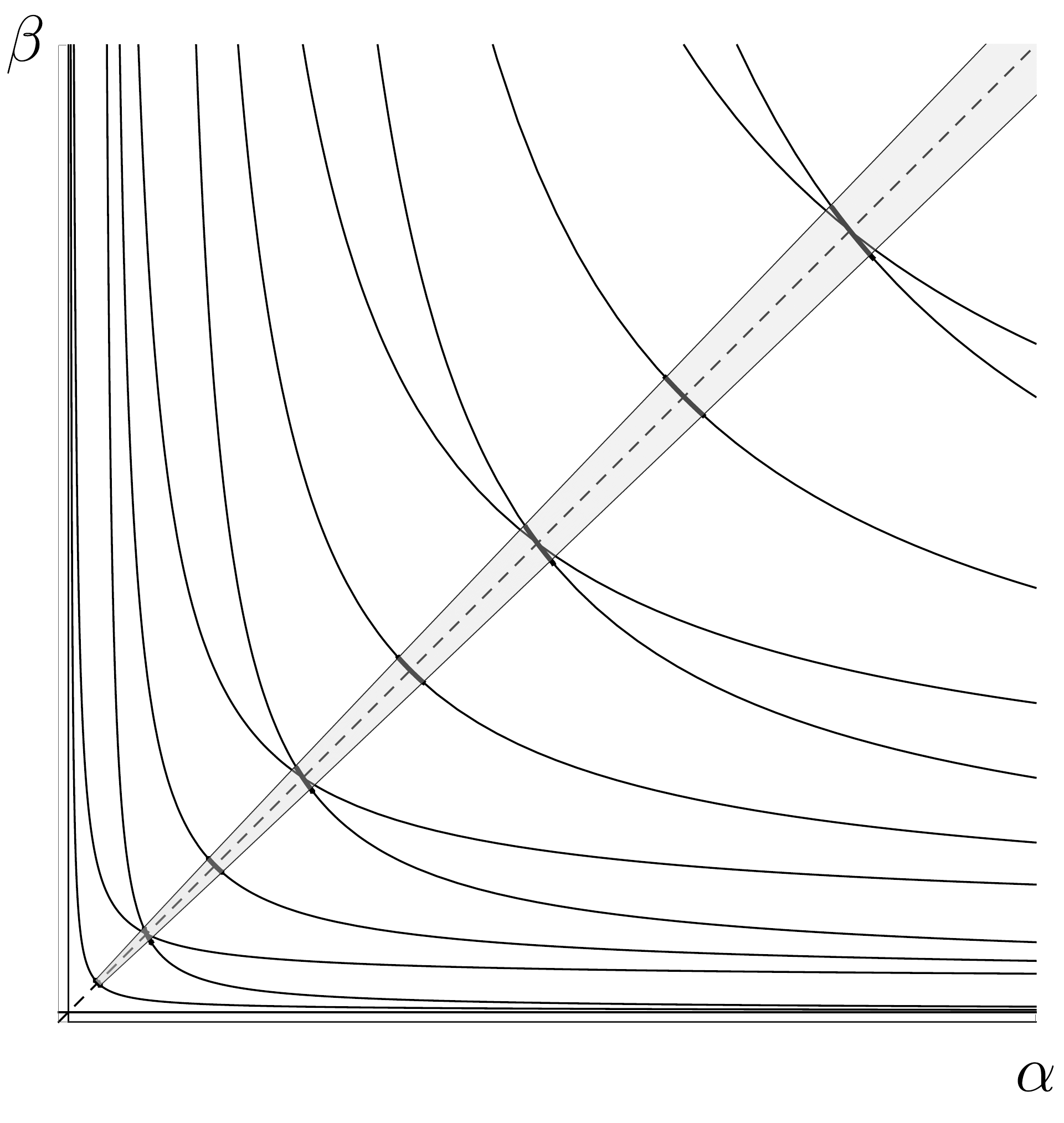}}}
	\qquad\qquad
	\subfloat[\centering $\epsilon = 0.5$]
	{{\includegraphics[width=0.32\linewidth]{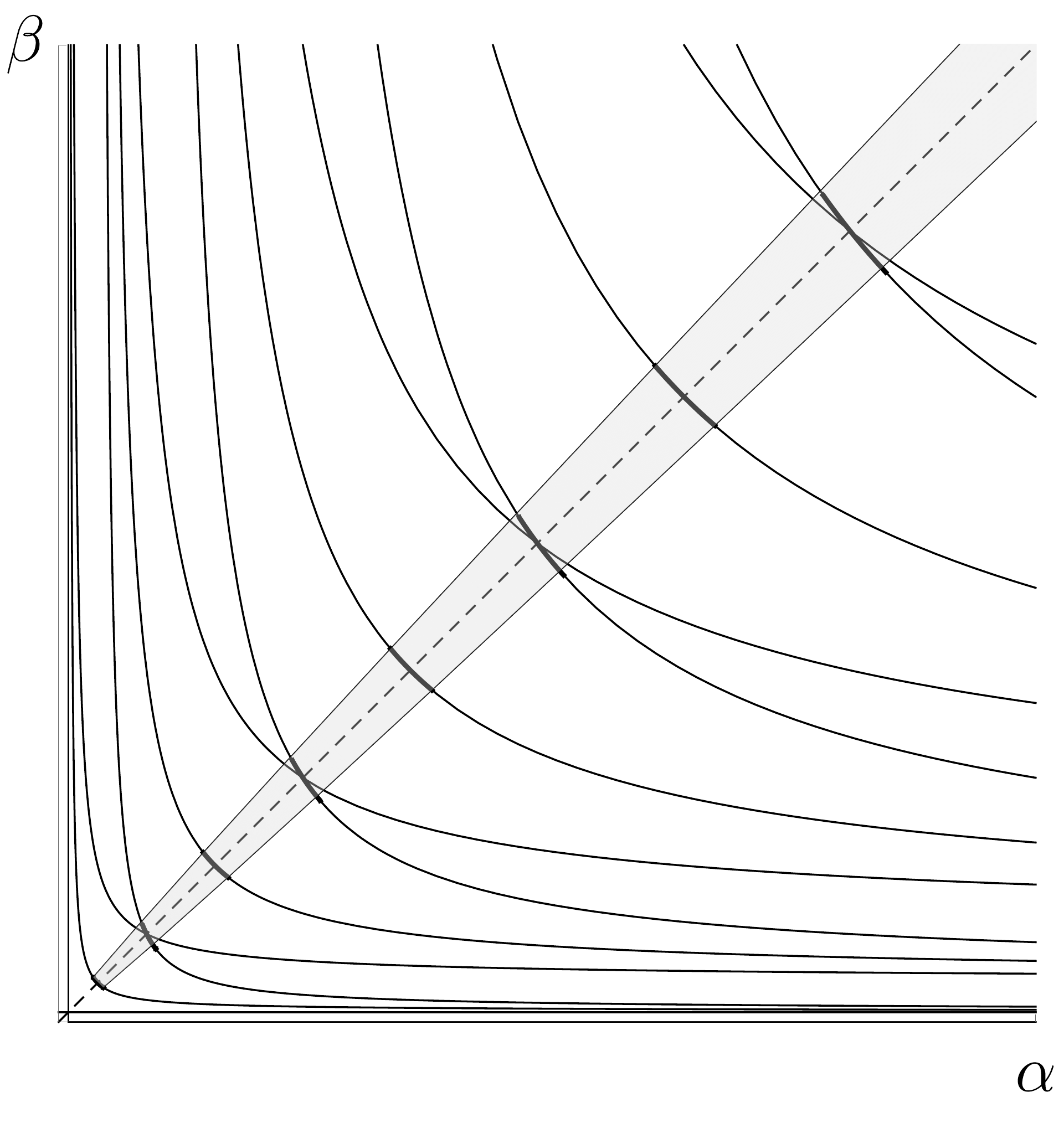}}}
	\caption{Shaded regions depict the result of Corollary \ref{cor:thm1} with $c_n = \frac{1}{46}\cdot\frac1{\zeta(1+\varepsilon)-1}$}
	\label{fig5}
\end{figure}

\begin{remark}
	The upper bounds \eqref{eq:bas_con2} and \eqref{eq:bas_con3} can be replaced by the following weaker ones which are independent from $n$:
	\begin{align}
	c_n&<\frac{1}{46}\cdot\frac1{\zeta(1+\varepsilon)-1} 
	\quad\quad \mbox{for odd }n\mbox{ with }\alpha(n)\geq\beta(n), \tag{\ref{eq:bas_con2}'} \\
	c_n&<\frac{1}{10}\cdot\frac1{\zeta(1+\varepsilon)-1} 
	\qquad \mbox{for odd }n\mbox{ with }\beta(n)>\alpha(n). \tag{\ref{eq:bas_con3}'}
	\end{align}
	On Figure \ref{fig5} we depict two regions of Fu\v{c}\'ik eigenvalues described by the inequality \eqref{eq:maxalpha} with the uniform constants $c_n$ given by (\ref{eq:bas_con2}').
\end{remark}

\medskip
Clearly, the general result given by Theorem \ref{th:main} also covers such Fu\v{c}\'ik systems $F_{\alpha,\beta}$ in which only  \textit{some} Fu\v{c}\'ik eigenfunctions differ from the sine functions. 
However, in the particular case when $f^n_{\alpha,\beta} = \phi_n$ for all odd $n$, we can obtain the following stronger result.
\begin{theorem}\label{thm:even}
Let $F_{\alpha,\beta}$ be a Fu\v{c}\'ik system with $\alpha(n)=\beta(n)=n^2$ for every odd $n$.
If the mappings $\alpha$ and $\beta$ satisfy
\begin{equation}\label{eq:thm:even}
\sum_{n=2}^\infty 
\left(\frac{\max(\sqrt{\alpha(n)},\sqrt{\beta(n)})}{n}-1\right)^2 
<
\infty,
\end{equation}
then $F_{\alpha,\beta}$ is a Riesz basis in $L^2(0,\pi)$.
\end{theorem}

Theorem \ref{thm:even} will be proven in Section \ref{sec:even} by showing that the corresponding Fu\v{c}\'ik system $F_{\alpha,\beta}$ is $\omega$-linearly independent and it is quadratically near to the system of sine functions $\{\phi_n\}$.

\begin{corollary}\label{cor:even}
	Let $F_{\alpha,\beta}$ be a Fu\v{c}\'ik system with $\alpha(n)=\beta(n)=n^2$ for every odd $n$.
	Let $\varepsilon,c>0$ be fixed. 
	If the mappings $\alpha$ and $\beta$ satisfy
	\begin{equation}\label{eq:cor-even}
		\max\left(\sqrt{\alpha(n)},\sqrt{\beta(n)}\right)\leq n
		+ \sqrt{c} \, n^{(1-\epsilon)/2}
	\end{equation}
	for any even $n$, 
	then $F_{\alpha,\beta}$ is a Riesz basis in $L^2(0,\pi)$.
\end{corollary}

\begin{figure}[h!]
	\centering
	\subfloat[\centering $\epsilon = 0.1$, $c=0.25$; $\gamma=5.6$]
	{{\includegraphics[width=0.32\linewidth]{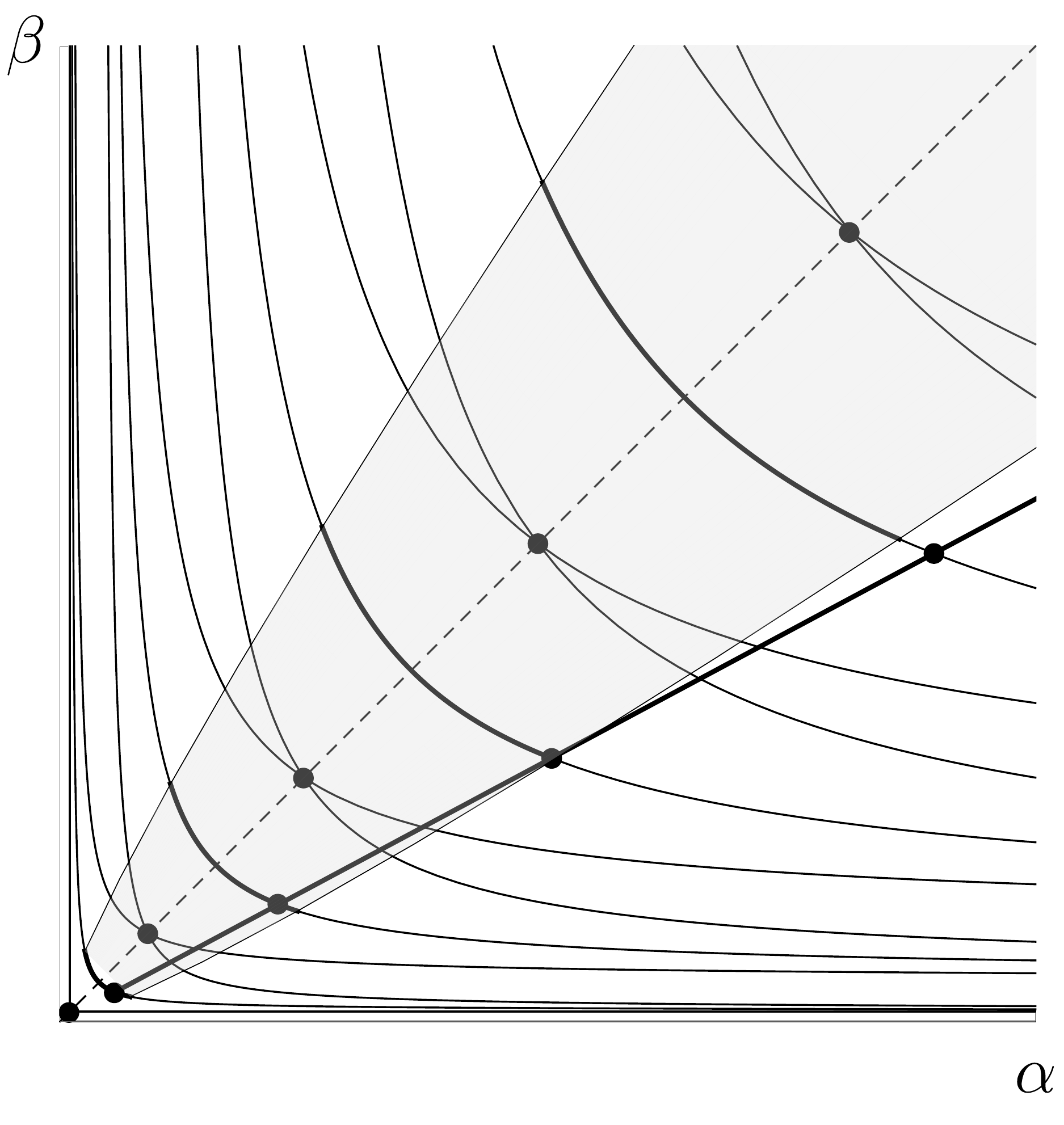}}}
	\qquad\qquad
	\subfloat[\centering $\epsilon = 0.5$, $c=0.4$; $\gamma=5.6$]
	{{\includegraphics[width=0.32\linewidth]{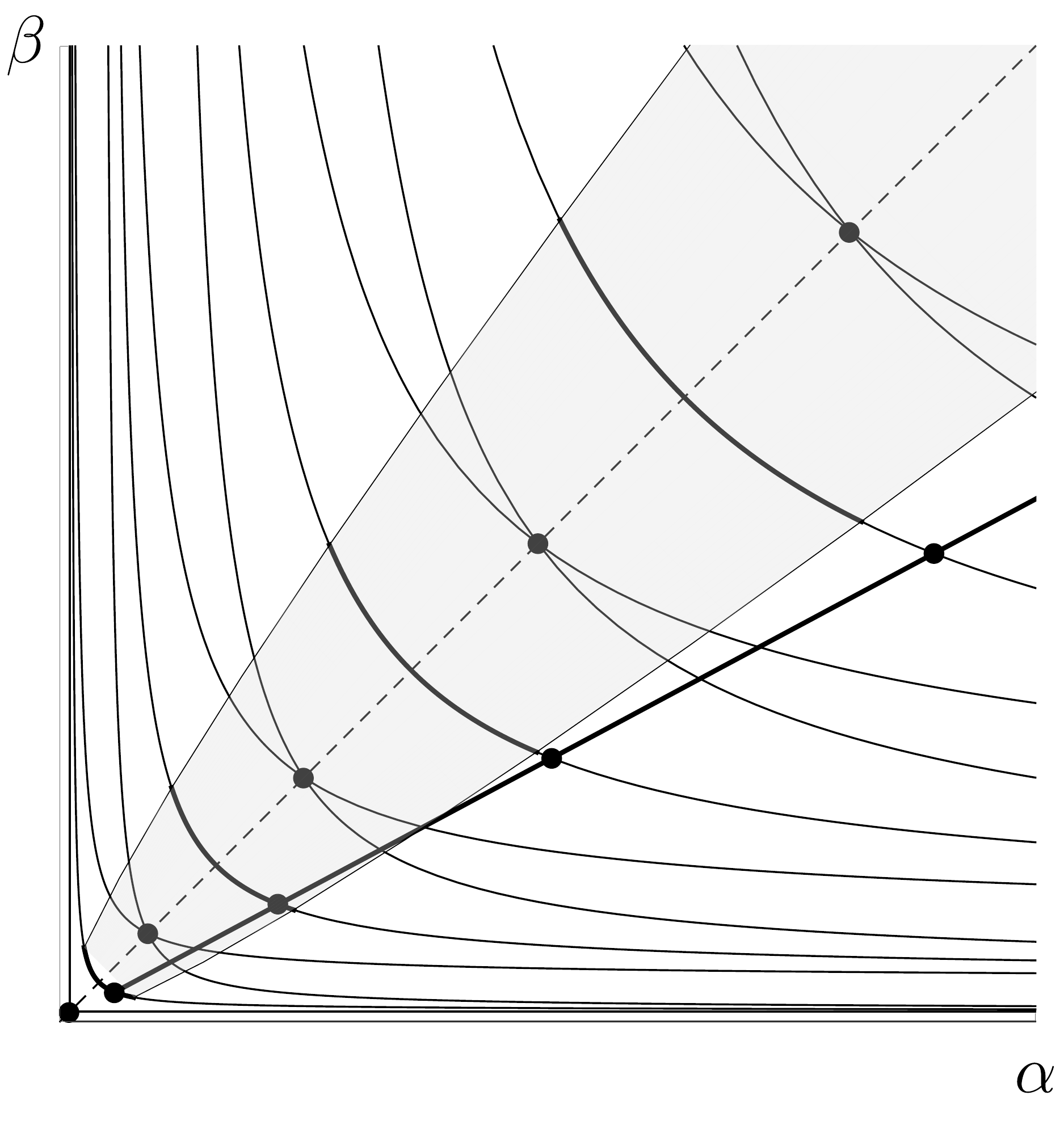}}}
	\caption{Comparison of the results of Corollary \ref{cor:even} and Theorem \ref{th:main2}}
	\label{fig4}
\end{figure}

\medskip
Let us now discuss the basisness of Fu\v{c}\'ik systems by means of the Paley-Wiener nearness to the system of sine functions.
We consider a Fu\v{c}\'ik system $F_{\alpha,\beta}$ for which the points $(\alpha(n),\beta(n))$ for even $n$ are on a line through the origin and, as in Theorem \ref{thm:even}, the normalized Fu\v{c}\'ik eigenfunctions $f^n_{\alpha,\beta}$ are just $\phi_n$ for odd $n$. 
We apply the method of separation of variables from \cite{duff} to this specific Fu\v{c}\'ik system to obtain mappings $\alpha$ and $\beta$ with better asymptotics as $n \to \infty$ than in Theorem \ref{thm:even}, see Figure \ref{fig4}.

\begin{theorem}\label{th:main2}
Let $F_{\alpha,\beta}$ be a Fu\v{c}\'ik system with $\alpha(n)=\beta(n)=n^2$ for every odd $n$ and
\begin{equation}\label{eq:DF1}
\alpha(n) = \frac{n^2 \gamma}4, \quad
\beta(n) = \frac{n^2 \gamma}{(2\sqrt{\gamma}-2)^2}
\end{equation}
for every even $n$, where $\gamma \in [4,5.682]$ is an arbitrary fixed constant. 
Then $F_{\alpha,\beta}$ is a Riesz basis in $L^2(0,\pi)$.
\end{theorem}

The choice of the Fu\v{c}\'ik eigenfunctions for even $n$ in Theorem \ref{th:main2} guarantees that these functions form a dilated system in the sense that
$$
f^n_{\alpha(n),\beta(n)}(x)
=
f^2_{\alpha(2),\beta(2)}\left(\frac{nx}2\right) \quad \text{for any even } n.
$$
This property is important in the proof of Theorem \ref{th:main2} which we give in Section \ref{sec:thm2}.

The formulas \eqref{eq:DF1} in Theorem \ref{th:main2} obviously guarantee that $\alpha(n)>\beta(n)$ for even $n$.
Moreover, the points $(\alpha(n),\beta(n))$ defined by \eqref{eq:DF1} are on the line 
\begin{equation}\label{eq:line}
\beta = \frac{4 \alpha}{(2\sqrt{\gamma}-2)^2}.
\end{equation}
We observe that \eqref{eq:DF1} for the mapping $\alpha$ can be written as
$$
\sqrt{\alpha(n)}= n+\left(\frac{\sqrt{\gamma}}2-1\right)n
\quad \text{for any even } n,
$$
which provides a better asymptotic than \eqref{eq:cor-even} of Corollary \ref{cor:even}.
Notice also that Theorem \ref{th:main2}
remains valid if we exchange $\alpha$ and $\beta$ due to the symmetry properties of the Fu\v{c}\'ik curves $\Gamma_{n}$ with even $n$.

Let us emphasize that constants which appear in the function $C_n$ in Theorem \ref{th:main} and the admissible range for the constant $\gamma$ in Theorem \ref{th:main2} are non-optimal since they follow from a series of estimates convenient for simplification of the derived expressions.
Thus, we anticipate that these constants might be substantially improved.
On the other hand, the asymptotic growth rate of $C_n(\alpha(n),\beta(n))$ is expected to be sharp for the quadratic nearness considerations.

\section{Norms -- distances -- scalar products}\label{sec:norms}

In this section, we derive explicit expressions for the $L^2$-norms of the normalized Fu\v{c}\'ik eigenfunctions $f^n_{\alpha,\beta}$, the distances $\|f^n_{\alpha,\beta}-\phi_n\|^2$ which are important for the proof of Theorem \ref{th:main}, and some scalar products $\langle f^n_{\alpha,\beta}, \phi_m \rangle$ which will be used in the proof of Theorem \ref{thm:even}.
We write these formulas for the case $\alpha\geq n^2\geq\beta$ in dependence of $\alpha$ and $n$, and for the case $\beta>n^2>\alpha$ in dependence of $\beta$ and $n$. 
We will thoroughly treat only the first case, and omit details for the second case to shorten the exposition.

Recall the following notations from Section \ref{sec:fucik}:
$$
l_1 = \frac{\pi}{\sqrt{\alpha}},
\quad
l_2 = \frac{\pi}{\sqrt{\beta}},
\quad
l = l_1 + l_2.
$$
Moreover, we have
\begin{equation*}
	\alpha = \frac{n^2\beta}{(2\sqrt{\beta}-n)^2} \quad\mbox{and}\quad \beta = \frac{n^2\alpha}{(2\sqrt{\alpha}-n)^2}
\end{equation*}
for all points $(\alpha,\beta)\in \Gamma_n$ with even $n$, and
\begin{equation*}
	\alpha = \frac{(n+1)^2\beta}{(2\sqrt{\beta}-(n-1))^2} \quad\mbox{and}\quad \beta = \frac{(n-1)^2\alpha}{(2\sqrt{\alpha}-(n+1))^2}
\end{equation*}
for all points $(\alpha,\beta)\in \Gamma_n$ with odd $n$.

\subsection{The case \texorpdfstring{$\alpha\geq n^2\geq\beta$}{a>n2>b}}\label{sec:norms:a>b}

We begin with the derivation of the norms of the normalized Fu\v{c}\'ik eigenfunctions $f^n_{\alpha,\beta}$. 
We obtain
\begin{align*}
	\|f^n_{\alpha,\beta}\|^2 &= 
	\frac{n}{2} \int_0^{l_1}\frac{\beta}{\alpha}\sin^2(\sqrt{\alpha}x)\,\mathrm{d}x 
	+ 
	\frac{n}{2}\int_0^{l_2}\sin^2(\sqrt{\beta}x)\,\mathrm{d}x
	= 	\frac{n}{4}\frac{\beta}{\alpha}\frac{\pi}{\sqrt{\alpha}} 
	+
	\frac{n}{4}\frac{\pi}{\sqrt{\beta}} \\
	&= 
	\frac{\pi n^3}{4\sqrt{\alpha}(2\sqrt{\alpha}-n)^2}
	+
	\frac{1}{2}\left(\pi-\frac{n}{2}\frac{\pi}{\sqrt{\alpha}}\right)
	= 
	\frac{\pi}{2} - \frac{\pi n(\sqrt{\alpha}-n)}{(2\sqrt{\alpha}-n)^2}
\end{align*}
for even $n$. 
In a similar way, we have
\begin{equation*}
	\|f^n_{\alpha,\beta}\|^2 = \frac{n+1}4\frac{\beta}{\alpha}\frac{\pi}{\sqrt{\alpha}} + \frac{n-1}4\frac{\pi}{\sqrt{\beta}}
	= \frac{\pi}2-\frac{\pi(n+1)(\sqrt{\alpha}-1)(\sqrt{\alpha}-n)}{\sqrt{\alpha}(2\sqrt{\alpha}-(n+1))^2}
\end{equation*}
for odd $n$.
Notice that all the sine functions have the same norm 
$\|\phi_n\|=\sqrt{\frac{\pi}{2}}$, $n\in\N$.

\medskip
Let us now derive expressions for the distances $\|f^n_{\alpha,\beta}-\phi_n\|^2$.
We start with the case $n=2$.
Using the formulas from Appendix \ref{sec:appendix}, we get
\begin{align}
	\notag
	\|
	f^2_{\alpha,\beta}-\phi_2
	\|^2 
	&=\int_0^{l_1}
	\left(
	\frac{\sqrt{\beta}}{\sqrt{\alpha}}\sin(\sqrt{\alpha}\,x)
	-
	\sin(2x)
	\right)^2
	\mathrm{d}x
	+ \int_{l_1}^\pi
	\left(
	\sin(\sqrt{\beta}(x-l_1)) + \sin(2x)
	\right)^2
	\mathrm{d}x
	\\
	\notag
	&= \frac{\pi}2 +\frac12\frac{\beta}{\alpha}l_1+\frac12l_2 + 2\frac{\sqrt{\beta}}{4-\alpha}\sin(2l_1) - 2\frac{\sqrt{\beta}}{4-\beta}\sin(2l_1) \\
	\notag
	&= \frac{\pi}2 +\frac12\frac{\beta}{\alpha}\frac{\pi}{\sqrt{\alpha}} + \frac12\frac{\pi}{\sqrt{\beta}} + 2\sqrt{\beta}\frac{\alpha-\beta}{(4-\alpha)(4-\beta)}\sin\left(\frac{2\pi}{\sqrt{\alpha}}\right)
	\\
	\label{eq:dist0}
	&=
	\pi - \pi\frac{2(\sqrt{\alpha}-2)}{(2\sqrt{\alpha}-2)^2} - \frac{4\alpha^2}{(2\sqrt{\alpha}-2)(3\sqrt{\alpha}-2)(\sqrt{\alpha}+2)}\frac{\sin\left(\frac{2\pi}{\sqrt{\alpha}}\right)}{\sqrt{\alpha}-2}.
\end{align}
Now, for a general even $n$ we have
\begin{align*}
	\|f^n_{\alpha,\beta}-\phi_n\|^2 
	&= 
	 \sum_{i=0}^{\frac{n}{2}-1}\int_{il}^{(i+1)l}(f^n_{\alpha,\beta}(x)-\sin(nx))^2\,\mathrm{d}x
	= \frac{n}{2}\int_0^l(f^n_{\alpha,\beta}(x)-\sin(nx))^2\,\mathrm{d}x \\
	&= \int_0^{\pi}(f^2_{4\alpha/n^2,4\beta/n^2}(y)-\sin(2y))^2\,\mathrm{d}y = \|f^2_{4\alpha/n^2,4\beta/n^2}-\phi_2\|^2,
\end{align*}
where we used that $f_{\alpha,\beta}(x)-\sin(nx)$ realizes the same values on each interval $(il,(i+1)l)$ for every $i\in\mathbb{N}$.
Therefore, we deduce from \eqref{eq:dist0} that
\begin{equation}\label{eq:dist-fn-even-a>b}
	\|f^n_{\alpha,\beta}-\phi_n\|^2 
	=
	\pi - \pi\frac{n(\sqrt{\alpha}-n)}{(2\sqrt{\alpha}-n)^2} - \frac{4\alpha^2}{(2\sqrt{\alpha}-n)(3\sqrt{\alpha}-n)(\sqrt{\alpha}+n)}\frac{\sin\left(\frac{n\pi}{\sqrt{\alpha}}\right)}{\sqrt{\alpha}-n}
\end{equation}
holds for any even $n$.

\medskip
Now we consider the case of odd $n\geq3$. 
This case requires more extensive calculations, and we put the derivation of the following formula to Appendix \ref{sec:appendix}.
We have
\begin{align}
\begin{split}\label{eq:dist-fn-odd-a>b}
\|f^n_{\alpha,\beta}-\phi_n\|^2
&= 
\pi - \pi\frac{(n+1)(\sqrt{\alpha}-1)}{\sqrt{\alpha}(2\sqrt{\alpha}-(n+1))^2}(\sqrt{\alpha}-n) \\
&- \frac{16(n-1)\sqrt{\alpha}^3}{(2\sqrt{\alpha}-(n+1))}\frac{(\sqrt{\alpha}-1)}{(n+\sqrt{\alpha})(n+1)((3n-1)\sqrt{\alpha}-n(n+1))}\times\\
&\times \frac{\cos\left(\frac{\pi}2\frac{n}{\sqrt{\alpha}}\right)\cos\left({\frac{\pi}{2}\frac{n^2+n-2\sqrt{\alpha}}{(n-1)\sqrt{\alpha}}}\right)}{(\sqrt{\alpha}-n)\sin\left({\pi\frac{\sqrt{\alpha}-n}{(n-1)\sqrt{\alpha}}}\right)}.
\end{split}
\end{align}
for any odd $n \geq 3$. 
Notice that the assumption $\alpha \geq n^2$ guarantees that each multiplier on the right-hand side of \eqref{eq:dist-fn-odd-a>b} is nonnegative.

\medskip
Finally, we derive some scalar products $\langle f^n_{\alpha,\beta}, \phi_m \rangle$ needed for the proof of Theorem \ref{thm:even}. For the special case $n=m \geq 2$, we can express the scalar product in terms of the formulas above as
\begin{equation}\label{eq:scalar}
\langle f^n_{\alpha,\beta},\phi_n\rangle = \frac12\left(\|f^n_{\alpha,\beta}\|^2+\left\|\phi_n\right\|^2
-
\|f^n_{\alpha,\beta}-\phi_n\|^2\right),
\end{equation}
and then easily obtain
\begin{equation}\label{eq:scalar-even1}
\langle f^n_{\alpha,\beta},\phi_n\rangle = \frac{2\alpha^2}{(2\sqrt{\alpha}-n)(3\sqrt{\alpha}-n)(\sqrt{\alpha}+n)}\frac{\sin\left(\frac{n\pi}{\sqrt{\alpha}}\right)}{\sqrt{\alpha}-n}
\end{equation}
for even $n$, 
and
\begin{align*}
\langle f^n_{\alpha,\beta},\phi_n\rangle &= \frac{8(n-1)\sqrt{\alpha}^3}{(2\sqrt{\alpha}-(n+1))}\frac{(\sqrt{\alpha}-1)}{(n+\sqrt{\alpha})(n+1)((3n-1)\sqrt{\alpha}-n(n+1))}\times\\
&\times \frac{\cos\left(\frac{\pi}2\frac{n}{\sqrt{\alpha}}\right)\cos\left({\frac{\pi}{2}\frac{n^2+n-2\sqrt{\alpha}}{(n-1)\sqrt{\alpha}}}\right)}{(\sqrt{\alpha}-n)\sin\left({\pi\frac{\sqrt{\alpha}-n}{(n-1)\sqrt{\alpha}}}\right)}
\end{align*}
for odd $n \geq 3$.

If $n \neq m$, then the scalar product $\langle f^n_{\alpha,\beta}, \phi_m \rangle$ vanishes for some combinations of $n$ and $m$. 
In particular, we have $\langle f^n_{\alpha,\beta}, \phi_m \rangle=0$ for odd $n$ and even $m$ by a simple symmetry argument. 
When both $n$ and $m$ are even with $n>m$, the scalar product also vanishes. 
Indeed, we get
\begin{equation}\label{eq:speon}
\langle f^n_{\alpha,\beta},\phi_m\rangle =
\frac{(\beta-\alpha)\sqrt{\beta}}{(\alpha-m^2)(\beta-m^2)}
\frac{\left(\sin\left(\left(\frac{n}{2}-1\right)\frac{m}{2}l\right)+\sin\left(\left(\frac{n}{2}-1\right)\frac{m}{2}l+ml_1\right)\right)\sin\left(\frac{mn}{4}l\right)}{\sin\left(\frac{m}{2}l\right)}
=0
\end{equation}
due to $\frac{n}2l=\pi$. 
Details on the derivation of this formula are given in Appendix \ref{sec:appendix}.

\subsection{The case \texorpdfstring{$\beta> n^2>\alpha$}{b>n2>a}}\label{sec:norms:b>a}
For the norm of the normalized Fu\v{c}\'ik eigenfunctions $f^n_{\alpha,\beta}$, we have
\begin{align*}\label{eq:norm-fn-even-b>a}
\begin{split}
	\|f^n_{\alpha,\beta}\|^2 &= \frac{n}{2}\int_0^{l_1}\sin^2(\sqrt{\alpha}x)\,\mathrm{d}x 
	+ 
	\frac{n}{2}\int_0^{l_2}\frac{\alpha}{\beta}\sin^2(\sqrt{\beta}x)\,\mathrm{d}x
	\\
	&= \frac{n}{4}\frac{\pi}{\sqrt{\alpha}} 
	+ 
	\frac{n}{4}\frac{\alpha}{\beta}\frac{\pi}{\sqrt{\beta}} = \frac{\pi}{2} - \frac{\pi n(\sqrt{\beta}-n)}{(2\sqrt{\beta}-n)^2}
\end{split}
\end{align*}
for even $n$, and 
\begin{equation*}\label{eq:norm-fn-odd-b>a}
	\|f^n_{\alpha,\beta}\|^2 = \frac{n+1}4\frac{\pi}{\sqrt{\alpha}} + \frac{n-1}4\frac{\alpha}{\beta}\frac{\pi}{\sqrt{\beta}}
	= \frac{\pi}2-\frac{\pi(n-1)(\sqrt{\beta}+1)(\sqrt{\beta}-n)}{\sqrt{\beta}(2\sqrt{\beta}-(n-1))^2}
\end{equation*}
for odd $n$.

\medskip
Now we derive expressions for the distances $\|f^n_{\alpha,\beta}-\phi_n\|^2$.
As in Section \ref{sec:norms:a>b}, 
we start with the case $n=2$:
\begin{align*}
	\|f^2_{\alpha,\beta}-\phi_2\|^2 &= \int_0^{l_1}(\sin(\sqrt{\alpha}\,x)-\sin(2x))^2\,\mathrm{d}x 
	+
	\int_{l_1}^\pi
	\left(
	\frac{\sqrt{\alpha}}{\sqrt{\beta}}\sin(\sqrt{\beta}(x-l_1))
	+
	\sin(2x)
	\right)^2\mathrm{d}x
	\\
	&=
	\frac{\pi}2 +\frac12\frac{\pi}{\sqrt{\alpha}} + \frac12\frac{\alpha}{\beta}\frac{\pi}{\sqrt{\beta}} + 2\sqrt{\alpha}\frac{\alpha-\beta}{(4-\alpha)(4-\beta)}\sin\left(\frac{2\pi}{\sqrt{\alpha}}\right)
	\\
	&=
	\pi - \pi\frac{2(\sqrt{\beta}-2)}{(2\sqrt{\beta}-2)^2} - \frac{4\beta^2}{(2\sqrt{\beta}-2)(3\sqrt{\beta}-2)(\sqrt{\beta}+2)}\frac{\sin\left(\frac{2\pi}{\sqrt{\beta}}\right)}{\sqrt{\beta}-2}.
\end{align*}
Therefore, we obtain
\begin{equation}\label{eq:dist-f-even-nb<a}
\|f^n_{\alpha,\beta}-\phi_n\|^2
=
\pi - \pi\frac{n(\sqrt{\beta}-n)}{(2\sqrt{\beta}-n)^2} - \frac{4\beta^2}{(2\sqrt{\beta}-n)(3\sqrt{\beta}-n)(\sqrt{\beta}+n)}\frac{\sin\left(\frac{n\pi}{\sqrt{\beta}}\right)}{\sqrt{\beta}-n}
\end{equation}
for even $n$.
Furthermore, we have
\begin{align}\label{eq:dist-f-odd-nb<a}
\begin{split}
	\|f^n_{\alpha,\beta}-\phi_n\|^2
	&= \pi - \pi\frac{(n-1)(\sqrt{\beta}+1)}{\sqrt{\beta}(2\sqrt{\beta}-(n-1))^2}(\sqrt{\beta}-n)\\
	&- \frac{16(n+1)\sqrt{\beta}^3}{(2\sqrt{\beta}-(n-1))}\frac{(\sqrt{\beta}+1)}{(n+\sqrt{\beta})(n-1)((3n+1)\sqrt{\beta}-n(n-1))}\times\\
	&\times \frac{\cos\left(\frac{\pi}2\frac{n}{\sqrt{\beta}}\right)\cos\left(\frac{\pi}{2}\frac{(2\sqrt{\beta}+n^2-n)}{(n+1)\sqrt{\beta}}\right)}{(\sqrt{\beta}-n)\sin\left(\pi\frac{n(\sqrt{\beta}+1)}{(n+1)\sqrt{\beta}}\right)}
\end{split}
\end{align}
for odd $n\geq3$. 
As in \eqref{eq:dist-fn-odd-a>b}, each multiplier on the right-hand side of \eqref{eq:dist-f-odd-nb<a} is nonnegative in view of the assumption $\beta > n^2$.

Finally, for the derivation of the scalar product $\langle f^n_{\alpha,\beta}, \phi_n \rangle$, we use \eqref{eq:scalar} and get
\begin{equation}\label{eq:scalar-even2}
\langle f^n_{\alpha,\beta},\phi_n\rangle = \frac{2\beta^2}{(2\sqrt{\beta}-n)(3\sqrt{\beta}-n)(\sqrt{\beta}+n)}\frac{\sin\left(\frac{n\pi}{\sqrt{\beta}}\right)}{\sqrt{\beta}-n}
\end{equation}
for even $n$, and
\begin{align*}
\langle f^n_{\alpha,\beta},\phi_n\rangle
	&=
\frac{8(n+1)\sqrt{\beta}^3}{(2\sqrt{\beta}-(n-1))}\frac{(\sqrt{\beta}+1)}{(n+\sqrt{\beta})(n-1)((3n+1)\sqrt{\beta}-n(n-1))}\times\\
	&\times \frac{\cos\left(\frac{\pi}2\frac{n}{\sqrt{\beta}}\right)\cos\left(\frac{\pi}{2}\frac{(2\sqrt{\beta}+n^2-n)}{(n+1)\sqrt{\beta}}\right)}{(\sqrt{\beta}-n)\sin\left(\pi\frac{n(\sqrt{\beta}+1)}{(n+1)\sqrt{\beta}}\right)}
\end{align*}
for odd $n \geq 3$.
The scalar product $\langle f^n_{\alpha,\beta}, \phi_m \rangle$ vanishes for even $m$ provided that either $n$ is odd, or $n$ is even with $n>m$, as in the case $\alpha \geq n^2 \geq \beta$ in Section \ref{sec:norms:a>b}.

\section{Asymptotics of distances}\label{sec:asympotics}

In this section, we establish upper bounds on the distances $\|f^n_{\alpha,\beta}-\phi_n\|^2$ which allows to describe their asymptotic behaviour for $n \to \infty$, as well as
for $\sqrt{\alpha} \to n$ or $\sqrt{\beta} \to n$.

\subsection{The case \texorpdfstring{$\alpha\geq n^2\geq\beta$}{a>n2>b}}\label{sec:asympotics:a>b}

We begin with estimating the distance $\|f^n_{\alpha,\beta}-\phi_n\|^2$ given by \eqref{eq:dist-fn-even-a>b} for even $n$. 
Using the lower bound
\begin{equation}\label{eq:sin-upper}
\sin x \geq x - \frac{x^3}{6} = \frac{x}{6}(\sqrt{6}-x)(\sqrt{6}+x), \quad 
x \geq 0,
\end{equation}
and the assumption $\alpha \geq n^2$, 
we see that
\begin{align*}
\sin\left(\frac{n\pi}{\sqrt{\alpha}}\right) 
=
\sin\left(\pi - \frac{\pi n}{\sqrt{\alpha}}\right)
\geq
\frac{\pi}{6 \sqrt{\alpha^3}}(\sqrt{\alpha}-n)
((\sqrt{6}-\pi)\sqrt{\alpha}+n\pi)((\sqrt{6}+\pi)\sqrt{\alpha}-n\pi).	
\end{align*}
Therefore, we get
\begin{align}
	\notag
\|f^n_{\alpha,\beta}-\phi_n\|^2
&\leq
\pi - \pi\frac{n(\sqrt{\alpha}-n)}{(2\sqrt{\alpha}-n)^2} 
- \frac{2\pi \sqrt{\alpha}	((\sqrt{6}-\pi)\sqrt{\alpha}+n\pi)((\sqrt{6}+\pi)\sqrt{\alpha}-n\pi)}{3(2\sqrt{\alpha}-n)(3\sqrt{\alpha}-n)(\sqrt{\alpha}+n)}
\\
	\notag
&
=
\frac{\pi}{3} \frac{4(3 + \pi^2)\alpha + \sqrt{\alpha} n (15 - 2 \pi^2) - 6 n^2}{(2\sqrt{\alpha}-n)^2(3\sqrt{\alpha}-n)(\sqrt{\alpha}+n)} (\sqrt{\alpha} - n)^2
\\
	\notag
&
\leq
\frac{\pi}{3} \frac{4(3 + \pi^2)\alpha - (2 \pi^2 - 9) n^2}{n^2(3\sqrt{\alpha}-n)(\sqrt{\alpha}+n)} (\sqrt{\alpha} - n)^2
\\
	\notag
&
=
\frac{4(3 + \pi^2)\pi}{9} \frac{\left(3\sqrt{\alpha} - \frac{3\sqrt{2 \pi^2 - 9}}{2\sqrt{3 + \pi^2}} n\right)\left(\sqrt{\alpha} + \frac{\sqrt{2 \pi^2 - 9}}{2\sqrt{3 + \pi^2}} n\right)}{n^2(3\sqrt{\alpha}-n)(\sqrt{\alpha}+n)} (\sqrt{\alpha} - n)^2\\
\label{eq:asymp:dist-fn-odd}
&
\leq
\frac{4(3 + \pi^2)\pi}{9}\frac{(\sqrt{\alpha}-n)^2}{n^2}
=
\frac{4(3 + \pi^2)\pi}{9}\left(\frac{\sqrt{\alpha}}{n} - 1\right)^2
\end{align}
for even $n$.

Let us now estimate $\|f^n_{\alpha,\beta}-\phi_n\|^2$ given by \eqref{eq:dist-fn-odd-a>b} for odd $n$. 
Recalling that all multipliers in \eqref{eq:dist-fn-odd-a>b} are nonnegative and applying the rough upper bound
$$
\sin\left(\pi\frac{\sqrt{\alpha}-n}{(n-1)\sqrt{\alpha}}\right)
\leq
\pi\frac{\sqrt{\alpha}-n}{(n-1)\sqrt{\alpha}},
$$
we get
\begin{align}
\begin{split}\label{eq:dist-est1}
\|f^n_{\alpha,\beta}-\phi_n\|^2 
&\leq \pi - \pi\frac{(n+1)(\sqrt{\alpha}-1)}{\sqrt{\alpha}(2\sqrt{\alpha}-(n+1))^2}(\sqrt{\alpha}-n) \\
&- \frac{16(n-1)^2\alpha^2}{\pi(2\sqrt{\alpha}-(n+1))}\frac{(\sqrt{\alpha}-1)}{(n+\sqrt{\alpha})(n+1)((3n-1)\sqrt{\alpha}-n(n+1))}\times\\
&\times \frac{\cos\left(\frac{\pi}2\frac{n}{\sqrt{\alpha}}\right)\cos\left(\frac{\pi}{2}\frac{n^2+n-2\sqrt{\alpha}}{(n-1)\sqrt{\alpha}}\right)}{(\sqrt{\alpha}-n)^2}.
\end{split}
\end{align}
Then, using \eqref{eq:sin-upper}, we have
$$
\cos(x) 
=
\sin\left(\frac{\pi}{2}-x\right)
\geq \frac{1}{6}\left(\frac{\pi}{2}-x\right)\left(\sqrt{6}-\frac{\pi}{2}+x\right)\left(\sqrt{6}+\frac{\pi}{2}-x\right),
\quad x \leq \frac{\pi}{2},
$$
and hence, by $\alpha>n^2$, we obtain
$$
\cos\left(\frac{\pi}2\frac{n}{\sqrt{\alpha}}\right) 
\geq
\frac{\pi(\sqrt{\alpha}-n)}{48\sqrt{\alpha^3}}((2\sqrt{6}-\pi)\sqrt{\alpha}+\pi n)((2\sqrt{6}+\pi)\sqrt{\alpha}-\pi n)
$$
and
\begin{align*}
\cos\left(\frac{\pi}{2}\frac{n^2+n-2\sqrt{\alpha}}{(n-1)\sqrt{\alpha}}\right)
\geq
\frac{\pi (n+1)(\sqrt{\alpha}-n)}{48 (n-1)^3 \sqrt{\alpha^3}}
&((2\sqrt{6}(n-1)-\pi (n+1))\sqrt{\alpha}+\pi n (n+1))
\\
&\times
((2\sqrt{6}(n-1) + \pi(n+1))\sqrt{\alpha}-\pi n (n+1)).
\end{align*}
Substituting these estimates into \eqref{eq:dist-est1}, we get
\begin{align}
\notag
\|f^n_{\alpha,\beta}-\phi_n\|^2 
&\leq 
\pi - \pi\frac{(n+1)(\sqrt{\alpha}-1)}{\sqrt{\alpha}(2\sqrt{\alpha}-(n+1))^2}(\sqrt{\alpha}-n) \\
\notag
&-\frac{\pi(\sqrt{\alpha}-1)((2\sqrt{6}-\pi)\sqrt{\alpha}+\pi n)((2\sqrt{6}+\pi)\sqrt{\alpha}-\pi n)}{144 (n-1) \alpha (2\sqrt{\alpha}-(n+1))(n+\sqrt{\alpha})((3n-1)\sqrt{\alpha}-n(n+1))}\times\\
\notag
&\times ((2\sqrt{6}(n-1)-\pi (n+1))\sqrt{\alpha}+\pi n (n+1))\\
\notag
&\times ((2\sqrt{6}(n-1) + \pi(n+1))\sqrt{\alpha}-\pi n (n+1))\\
\label{eq:asymp-fn-odd-1}
&\leq
\frac{4\pi n^2(n^2+1)}{(n-1)^4} \left(\frac{\sqrt{\alpha}}{n}-1\right)^2
\leq 
23 \pi \left(\frac{\sqrt{\alpha}}{n}-1\right)^2
\end{align}
for odd $n$.
Several intermediate estimates used to derive these upper bounds are given in Appendix \ref{sec:appendix}.

\subsection{The case \texorpdfstring{$\beta >  n^2 > \alpha$}{b>n2>a}}

Let us estimate $\|f^n_{\alpha,\beta}-\phi_n\|^2$ given by \eqref{eq:dist-f-even-nb<a} for even $n$. 
Noting that \eqref{eq:dist-f-even-nb<a} is the same formula as \eqref{eq:dist-fn-even-a>b} up to a replacement of $\alpha$ by $\beta$, we get 
\begin{equation}\label{eq:asymp:dist-fn-odd2}
\|f^n_{\alpha,\beta}-\phi_n\|^2
\leq
\frac{4(3 + \pi^2)\pi}{9}\left(\frac{\sqrt{\beta}}{n} - 1\right)^2
\end{equation}
for even $n$, 
as in \eqref{eq:asymp:dist-fn-odd}.

Now we provide an upper bound on $\|f^n_{\alpha,\beta}-\phi_n\|^2$ given by \eqref{eq:dist-f-odd-nb<a} for odd $n$. 
Using \eqref{eq:sin-upper} as in Section \ref{sec:asympotics:a>b} and recalling that $\beta>n^2$, we estimate the trigonometric terms in \eqref{eq:dist-f-odd-nb<a} as follows:
\begin{align*}
\cos\left(\frac{\pi}2\frac{n}{\sqrt{\beta}}\right) 
&\geq
\frac{\pi(\sqrt{\beta}-n)}{48\sqrt{\beta^3}}((2\sqrt{6}-\pi)\sqrt{\beta}+\pi n)((2\sqrt{6}+\pi)\sqrt{\beta}-\pi n),
\\
\cos\left(\frac{\pi}{2}\frac{(2\sqrt{\beta}+n^2-n)}{(n+1)\sqrt{\beta}}\right)
&\geq
\frac{\pi (n-1)(\sqrt{\beta}-n)}{48 (n+1)^3 \sqrt{\beta^3}}
((2\sqrt{6}(n+1)-\pi (n-1))\sqrt{\beta}+\pi n (n-1))
\\
&\times
((2\sqrt{6}(n+1) + \pi(n-1))\sqrt{\beta}-\pi n (n-1)),
\\
\sin\left(\pi\frac{n(\sqrt{\beta}+1)}{(n+1)\sqrt{\beta}}\right) 
&=
\sin\left(\pi - \pi\frac{n(\sqrt{\beta}+1)}{(n+1)\sqrt{\beta}}\right) 
=
\sin\left(\frac{\pi(\sqrt{\beta}-n)}{(n+1)\sqrt{\beta}}\right) 
\leq 
\frac{\pi(\sqrt{\beta}-n)}{(n+1)\sqrt{\beta}}.
\end{align*}
Substituting these estimates into \eqref{eq:dist-f-odd-nb<a}, we deduce that
\begin{align}
\notag
\|f^n_{\alpha,\beta}-\phi_n\|^2
&\leq 
\pi - \pi\frac{(n-1)(\sqrt{\beta}+1)}{\sqrt{\beta}(2\sqrt{\beta}-(n-1))^2}(\sqrt{\beta}-n)\\
\notag
&- \frac{\pi (\sqrt{\beta}+1) ((2\sqrt{6}-\pi)\sqrt{\beta}+\pi n)((2\sqrt{6}+\pi)\sqrt{\beta}-\pi n)}{144  (n+1)\beta (2\sqrt{\beta}-(n-1))(n+\sqrt{\beta})((3n+1)\sqrt{\beta}-n(n-1))}
\\
\notag
&\times
((2\sqrt{6}(n+1)-\pi (n-1))\sqrt{\beta}+\pi n (n-1))
\\
\notag
&\times
((2\sqrt{6}(n+1) + \pi(n-1))\sqrt{\beta}-\pi n (n-1))
\\
\label{eq:asymp:dist-fn-odd-b>a}
&\leq
\frac{5 \pi n^2(n^2+1)}{(n+1)^4}
\left(\frac{\sqrt{\beta}}{n}-1\right)^2
\leq
5 \pi \left(\frac{\sqrt{\beta}}{n}-1\right)^2
\end{align}
for odd $n$.
Several intermediate estimates in the derivation of these upper bounds can be found in Appendix \ref{sec:appendix}.

\section{Proof of Theorem \ref{thm:even}}\label{sec:even}

To prove Theorem \ref{thm:even}, we first establish the $\omega$-linear independence of Fu\v{c}\'ik systems $F_{\alpha,\beta}$ with $\alpha(n)=\beta(n)=n^2$ for every odd $n$.
As the second step, we use the bounds \eqref{eq:asymp:dist-fn-odd} and \eqref{eq:asymp:dist-fn-odd2} to deduce that the assumption \eqref{eq:thm:even} guarantees the quadratic nearness of $F_{\alpha,\beta}$ to the system of sine functions $\{\phi_n\}$.

Let $\{\eta_n\}$ be a sequence of scalars such that 
$$
\lim_{m \to \infty} 
\|\sum_{n=1}^{m}\eta_n f^n_{\alpha,\beta}\|=0.
$$
In order to show that $F_{\alpha,\beta}$ is $\omega$-linear independent, we have to prove $\eta_n=0$ for each $n$.
Fix any $k \in \mathbb{N}$ and consider the functions
$$
g_{m,k}(x) = \sin (kx) \sum_{n=1}^{m}\eta_n f^n_{\alpha,\beta}(x)
\quad 
\text{for}~ m \in \mathbb{N}.
$$
We see that
$$
\Big| 
\int_0^\pi g_{m,k}(x) \, \mathrm{d}x
\Big|
\leq
\int_0^\pi |\sin(kx)| \Big| \sum_{n=1}^{m}\eta_n f^n_{\alpha,\beta}(x)\Big|  \, \mathrm{d}x \leq 
\sqrt{\frac{\pi}{2}}
\, \|\sum_{n=1}^{m}\eta_n f^n_{\alpha,\beta}\| \to 0
\quad \text{as}~ m \to \infty,
$$
which yields 
\begin{equation}\label{eq:sumlimit}
\sum_{n=1}^{\infty}\eta_n 
\langle
f^n_{\alpha,\beta}, \phi_k
\rangle
=
\lim_{m \to \infty} 
\int_0^\pi g_{m,k}(x) \, \mathrm{d}x
=
0.
\end{equation}
Taking $k = 2$ and recalling that  $\langle f^n_{\alpha,\beta}, \phi_2 \rangle = 0$ for any $n \neq 2$ (see Section \ref{sec:norms}), we get $\eta_2 \langle f^2_{\alpha,\beta}, \phi_2 \rangle =0$. 
Since $\langle f^2_{\alpha,\beta}, \phi_2 \rangle \neq 0$ by \eqref{eq:scalar-even1} and \eqref{eq:scalar-even2}, we obtain $\eta_2=0$.
By an inductive argument, we derive in the same way $\eta_{2m}=0$ for every $m \geq 1$.
Let us remark that until now we did not use the special form of the Fu\v{c}\'ik system $F_{\alpha,\beta}$, namely,  $\alpha(n)=\beta(n)=n^2$ for all odd $n \geq 3$.
Since this assumption reads as $f^n_{\alpha,\beta}=\phi_n$ for any odd $n$, we further obtain from \eqref{eq:sumlimit} that $\eta_{2m-1}=0$ for every $m \geq 1$ by the orthogonality of $\{\phi_n\}$.
This concludes the $\omega$-linear independence of the Fu\v{c}\'ik system $F_{\alpha,\beta}$.

Using the uppers bounds \eqref{eq:asymp:dist-fn-odd} and \eqref{eq:asymp:dist-fn-odd2} on the distances $\|f^n_{\alpha,\beta}-\phi_n\|^2$, we get
$$
\sum_{n=1}^\infty 
\|f^n_{\alpha,\beta}-\phi_n\|^2
\leq 
\frac{4(3 + \pi^2)\pi}{9}
\sum_{n=2}^\infty 
\left(\frac{\max(\sqrt{\alpha(n)},\sqrt{\beta(n)})}{n}-1\right)^2.
$$
Since the right-hand side is bounded in view of the assumption \eqref{eq:thm:even}, the Fu\v{c}\'ik system $F_{\alpha,\beta}$ is quadratically near to the system of sine functions $\{\phi_n\}$.
Therefore, the rescaled system $\sqrt{2/\pi} \, F_{\alpha,\beta}$ is quadratically near to the complete orthonormal system $\{\sqrt{2/\pi} \, \phi_n\}$, and hence it is a Riesz basis in $L^2(0,\pi)$ by \cite[Theorem V-2.20]{kato}.
Clearly, $F_{\alpha,\beta}$ is also a Riesz basis.

\section{Proof of Theorem \ref{th:main2}}\label{sec:thm2}

In this section, we provide the proof of our third main result, 
Theorem \ref{th:main2}. 
For this purpose, we use the method of separation of variables from \cite{duff} to show that a specific Fu\v{c}\'ik system $F_{\alpha,\beta}$ is Paley-Wiener near to the system of sine functions $\{\phi_n\}$, see Definition \ref{def:FS}.
The classical result of Paley and Wiener \cite{paley} then yields the basisness of $F_{\alpha,\beta}$ in $L^2(0,\pi)$.

Recall that we choose $F_{\alpha,\beta}$ to satisfy $f^n_{\alpha(n),\beta(n)} = \phi_n$ for any odd $n$, and 
\begin{equation}
	\tag{\ref{eq:DF1}}
	\alpha(n) = \frac{n^2 \gamma}4, \quad
	\beta(n) = \frac{n^2 \gamma}{(2\sqrt{\gamma}-2)^2}
\end{equation}
for every even $n$, where $\gamma$ is a fixed constant  in the interval $[4,5.682]$.
Denoting $l_1(n) = \frac{\pi}{\sqrt{\alpha(n)}}$ and $l_2(n) = \frac{\pi}{\sqrt{\beta(n)}}$, we see that \eqref{eq:DF1} implies
\begin{equation*}\label{eq:l1nl123}
	\frac{l_1(2)}{l_1(n)}
	=
	\frac{n}{2}
	\quad \text{and} \quad
	\frac{l_2(2)}{l_2(n)}
	=
	\frac{n}{2}
\end{equation*}
for every even $n$.
It can be easily deduced from the piecewise definition \eqref{eq:fucikpiecewise1} that $f^n_{\alpha(n),\beta(n)}$ for even $n$ form a sequence of dilated functions in the sense that
\begin{equation}\label{eq:f-even1}
	f_{\alpha(n),\beta(n)}^n(x) = 
	f_{\alpha(2),\beta(2)}^2\left(\frac{nx}{2}\right)
	\equiv
	f_{\gamma,\gamma/(\sqrt{\gamma}-1)^2}^2\left(\frac{nx}{2}\right).
\end{equation}

We claim that the assumptions of \cite[Theorem D]{duff} (with $f_n = \phi_n$ and $g_n = f_{\alpha(n),\beta(n)}^n$) are satisfied for any $\gamma \in [4,5.682]$, 
namely, there exist a matrix of constants $\{C_{n,k}\}$ and a sequence of bounded linear operators $\{T_k\}$ such that each $f_{\alpha(n),\beta(n)}^n$ has the representation
\begin{equation}\label{eq:DE1}
	f_{\alpha(n),\beta(n)}^n(x) = \sin(nx) + \sum_{k=1}^{\infty}C_{n,k} T_k\sin(nx),
\end{equation}
where we have $|C_{n,k}| \leq c_k$ and $\|T_k\|_* \leq t_k$ with constants $c_k$ and $t_k$ satisfying $\sum_{k=1}^{\infty} c_k t_k < 1$.
Here, $\|\cdot\|_*$ is the operator norm.
If our claim is true, then the Fu\v{c}\'ik system $F_{\alpha,\beta}$ is a Riesz basis in $L^2(0,\pi)$.
We remark that although the system $\{\phi_n\}$ is not orthonormal, \cite[Theorem D]{duff} is applicable by simple rescaling arguments since all $\phi_n$ are of the same norm.

Let $n$ be even.
Then, in view of \eqref{eq:f-even1}, we have
\begin{equation}\label{eq:fourierfn}
	f_{\alpha(n),\beta(n)}^{n}(x) = f_{\alpha(2),\beta(2)}^2\left(\frac{nx}{2}\right) = \sum_{k=1}^{\infty}A_k\sin\left(\frac{knx}{2}\right),
\end{equation}
where $A_k$ are the coefficients of the odd Fourier expansion of $f_{\alpha(2),\beta(2)}^2$.
The comparison of \eqref{eq:DE1} and \eqref{eq:fourierfn} suggests to define bounded linear operators $T_k$, $k \in \mathbb{N}$, that satisfy the following property:	
\begin{equation}\label{eq:tksin}
T_k \sin(nx) = \sin\left(\frac{knx}{2}\right)
\quad \text{for every even}~ n.
\end{equation}
To this end, for any $g \in L^2(0,\pi)$ we define its  antiperiodic extension $g^*$ as
$$
g^*(x) = (-1)^\kappa g(x-\pi \kappa) \quad \text{for}~~ \pi \kappa \leq x \leq \pi (\kappa+1), \quad \kappa \in \mathbb{N}_0.
$$
In particular, we see that if $g(x) = \sin(x)$ for $x \in (0,\pi)$, then $g^*(x) = \sin(x)$ for any $x \in \mathbb{R}_+$.
Now we choose $T_k : L^2(0,\pi) \to L^2(0,\pi)$ as
\begin{equation}\label{eq:TK}
T_k g(x) = g^*\left(\frac{kx}{2}\right).
\end{equation}
Clearly, $T_2$ is just the identity operator and each $T_k$ satisfies \eqref{eq:tksin}.
Moreover, each $T_k$ is linear with the norm $\|T_k\|_* = 1$ for even $k$ and $\|T_k\|_* = \sqrt{1+1/k}$ for odd $k$, see Appendix \ref{sec:appendix2}.

Thus, in accordance with \eqref{eq:DE1} and \eqref{eq:fourierfn}, for even $n$ we set $C_{n,1}=A_1$, $C_{n,2}=A_2-1$, and $C_{n,k}=A_k$ for $k \geq 3$,
while for odd $n$ we simply choose $C_{n,k}=0$ for $k \in \mathbb{N}$.
For this choice of constants, we can set $c_1 = |A_1|$, $c_2=|A_2-1|$, and $c_k = |A_k|$ for $k \geq 3$.

Let us now estimate the constants $c_k$.
Notice that
\begin{align*}
	A_k = \frac2{\pi}\int_0^{\pi}f_{\alpha(2),\beta(2)}^2(x)\sin(kx)\,\mathrm{d}x  
	= \frac2{\pi}\frac{\gamma^2}{\sqrt{\gamma}-1}\frac{(2-\sqrt{\gamma})\sin\left(\frac{k\pi}{\sqrt{\gamma}}\right) }{(k^2-\gamma)(k^2(\sqrt{\gamma}-1)^2-\gamma)}.
\end{align*}
Thus, for $c_1$ we easily obtain
$$
c_1 = |A_1|
\leq
\frac2{\pi}\frac{\gamma^2(\sqrt{\gamma}-2)}{(\sqrt{\gamma}-1)^2(\sqrt{\gamma}+1)(2\sqrt{\gamma}-1)}.
$$
For $c_2$ we recall that $\gamma \geq 4$ and use the upper bound
\begin{equation}\label{eq:sinlower}
\sin\left(\frac{2\pi}{\sqrt{\gamma}}\right)=\sin\left(\pi-\frac{2\pi}{\sqrt{\gamma}}\right)=\sin\left(\frac{\pi(\sqrt{\gamma}-2)}{\sqrt{\gamma}}\right)\leq\frac{\pi(\sqrt{\gamma}-2)}{\sqrt{\gamma}}
\end{equation}
to deduce
\begin{align*}
	A_2-1 &= \frac2{\pi}\frac{\gamma^2}{\sqrt{\gamma}-1}\frac{\sin\left(\frac{2\pi}{\sqrt{\gamma}}\right)}{(\gamma-4)(3\sqrt{\gamma}-2)}-1 \leq \frac{2\sqrt{\gamma}^3}{(\sqrt{\gamma}-1)(\sqrt{\gamma}+2)(3\sqrt{\gamma}-2)}-1 \leq 0.
\end{align*}
Hence, using \eqref{eq:sin-upper} instead of \eqref{eq:sinlower}, we derive that
\begin{align*}
	c_2 
	= 1-A_2 
	&\leq 1-\frac{\sqrt{\gamma}((\sqrt{6}-\pi)\sqrt{\gamma}+2\pi)((\sqrt{6}+\pi)\sqrt{\gamma}-2\pi)}{3(\sqrt{\gamma}-1)(\sqrt{\gamma}+2)(3\sqrt{\gamma}-2)} \\
	&= \frac{((3+\pi^2)\gamma+(9-2\pi^2)\sqrt{\gamma}-6)(\sqrt{\gamma}-2)}
	{3(\sqrt{\gamma}-1)(\sqrt{\gamma}+2)(3\sqrt{\gamma}-2)}.
\end{align*}
Finally, for $c_k$ with $k \geq 3$ we get
\begin{align*}
c_k 
= |A_k|
&= \frac{2}{\pi}\frac{\gamma^2(\sqrt{\gamma}-2)}{\sqrt{\gamma}-1}
\frac{
	\big|
	\sin\left(\frac{k\pi}{\sqrt{\gamma}}\right)
	\big|
}{(k^2-\gamma)(k^2-\gamma+k^2\sqrt{\gamma}(\sqrt{\gamma}-2))} \leq \frac{2}{\pi}\frac{\gamma^2(\sqrt{\gamma}-2)}{\sqrt{\gamma}-1}\frac1{(k^2-\gamma)^2}.
\end{align*}

Using the estimates above, we deduce that
\begin{align*}
\sum_{k=1}^{\infty} c_k\|T_k\|_*  
&\leq
\sum_{k=1}^{4} c_k\|T_k\|_*  
+
\sqrt{\frac65} \cdot \frac{2}{\pi}\frac{\gamma^2(\sqrt{\gamma}-2)}{\sqrt{\gamma}-1} \,\sum_{k=5}^{\infty}\frac1{(k^2-9)^2}\\
&\leq \sqrt{2}\cdot\frac2{\pi}\frac{\gamma^2(\sqrt{\gamma}-2)}{(\sqrt{\gamma}-1)^2(\sqrt{\gamma}+1)(2\sqrt{\gamma}-1)} \\
&+ \frac{((3+\pi^2)\gamma+(9-2\pi^2)\sqrt{\gamma}-6)(\sqrt{\gamma}-2)}{3(\sqrt{\gamma}-1)(\sqrt{\gamma}+2)(3\sqrt{\gamma}-2)} \\
&+ 
\sqrt{\frac43} \cdot \frac{2}{\pi}\frac{\gamma^2(\sqrt{\gamma}-2)}{\sqrt{\gamma}-1}\frac1{(3^2-\gamma)^2} + \frac{2}{\pi}\frac{\gamma^2(\sqrt{\gamma}-2)}{\sqrt{\gamma}-1}\frac1{(4^2-\gamma)^2} \\
&+ 
\sqrt{\frac65} \cdot \frac{2}{\pi}\frac{\gamma^2(\sqrt{\gamma}-2)}{\sqrt{\gamma}-1}\left(\frac{\pi^2}{108}-\frac{536741}{6350400}\right)
=:
E(\gamma).
\end{align*}
By straightforward calculations, it is not hard to show that each summand in $E(\gamma)$ is strictly increasing with respect to $\gamma \geq 4$ and $E(4)=0$.
At the same time, we have $E(5.682) = 0.9992...$, which shows that $E(\gamma)<1$ for any $\gamma \in [4,5.682]$.
This completes the proof of Theorem \ref{th:main2}.

\section{Final remarks}\label{sec:final}

1. 
The quadratic nearness assumption \eqref{eq:sqna} used in the proof of Theorem \ref{th:main} can be weakened to the inequality
\begin{equation}\label{eq:assum_weak}
	\sum_{n=1}^{\infty}
	\left(
	\|f^n_{\alpha,\beta}-\phi_n\|^2
	-
	\frac{|\langle f^n_{\alpha,\beta}-\phi_n ,f^n_{\alpha,\beta}\rangle|^2}{\|f^n_{\alpha,\beta}\|^2}
	\right)
	<
	\frac{\pi}{2}
\end{equation}
and still guarantee that the Fu\v{c}\'ik system $F_{\alpha,\beta}$ is a Riesz basis in $L^2(0,\pi)$, see \cite[Theorem V-2.21]{kato}.
Noting that each summand in \eqref{eq:assum_weak} can be written as 
\begin{equation}\label{eq:kato2}
	\|f^n_{\alpha,\beta}-\phi_n\|^2-\frac{\left( \|f^n_{\alpha,\beta}\|^2-\|\phi_n\|^2+\|f^n_{\alpha,\beta}-\phi_n\|^2\right)^2}{4\|f^n_{\alpha,\beta}\|^2},
\end{equation}
one can apply the formulas from Section \ref{sec:norms} to derive the explicit expression for \eqref{eq:kato2} and estimate it from above in the same way as in Section \ref{sec:asympotics}. 
However, this does not improve the asymptotic behaviour of the function $C_n(\alpha(n),\beta(n))$ as $n \to \infty$ in Theorem \ref{th:main}, but only slightly improves the constants.

2.
We anticipate that the $\omega$-linear independence is satisfied for general Fu\v{c}\'ik systems. 
Provided this claim is true, 
the assumptions of Theorem \ref{th:main} can be significantly weakened in the sense that the sum in \eqref{eq:sumBn} only has to converge. 

3. 
In the proof of Theorem \ref{th:main2} we used the fact that normalized Fu\v{c}\'ik eigenfunctions $f^n_{\alpha,\beta}$, with $(\alpha(n),\beta(n))$ for even $n$ on a straight line through the origin, form a dilated sequence in the sense of \eqref{eq:f-even1}. 
It is not hard to show that the normalized Fu\v{c}\'ik eigenfunctions for odd $n$ with $(\alpha(n),\beta(n))$ on the same line satisfy
$$
f^n_{\alpha(n),\beta(n)}(x)
=
f^2_{\alpha(2),\beta(2)}
\left(
\left(\frac{n-1}{2}+ \frac{1}{\sqrt{\alpha(2)}}
\right)x
\right).
$$
Nevertheless, this property is less suitable for application of \cite[Theorem D]{duff} since the operators $T_k$ in \eqref{eq:DE1} have to be independent of $n$.

4.
Theorem \ref{th:main2} might suggest that the basisness of a Fu\v{c}\'ik system $F_{\alpha,\beta}$ also holds when each point $(\alpha(n), \beta(n)) \in \Gamma_n$, $n \geq 2$, belongs to the angular sector \textit{in between} the line \eqref{eq:line} and its reflection with respect to the main diagonal $\alpha=\beta$.

5.
Recall that the constants in the function $C_n$ in Theorem \ref{th:main}, as well as the upper bound for $\gamma$ in Theorem \ref{th:main2}, are not optimal due to the employed methods and the estimation procedure.

\appendix
\section{}\label{sec:appendix}

We start by providing several useful formulas.
For the derivation of expressions in Section \ref{sec:norms}, we need to calculate several integrals of the general form
\begin{equation*}\label{eq:int}
	\int
	\left(c\sin(\sqrt{\delta}(x-x_0))
	\pm\sin(nx)
	\right)^2
	\mathrm{d}x
\end{equation*}
with constants $c\in\left\{1,\sqrt{\alpha}/\sqrt{\beta},\sqrt{\beta}/\sqrt{\alpha}\right\}$, $\delta\in\{\alpha,\beta\}$, and certain $x_0\in[0,\pi]$. 
The antiderivative of this integral can be expressed through the following formulas:
\begin{align*}
	\int\sin^2(nx)\,\mathrm{d}x &= \frac12\left(x-\frac1{n}\sin(nx)\cos(nx)\right) + C, \\
	\int\sin^2(\sqrt{\delta}(x-x_0))\,\mathrm{d}x &= \frac12\left(x-x_0-\frac1{\sqrt{\delta}}\sin(\sqrt{\delta}(x-x_0))\cos(\sqrt{\delta}(x-x_0))\right)+C, \\
	\int \sin(\sqrt{\delta}(x-x_0))\sin(nx) \,\mathrm{d}x 
	&= \frac{\sqrt{\delta}}{n^2-\delta}
	\cos(\sqrt{\delta}(x-x_0))\sin(nx) \\
	&- \frac{n}{n^2-\delta}
	\sin(\sqrt{\delta}(x-x_0))\cos(nx)
	+C.
\end{align*}
We observe that $\sin(\sqrt{\delta}(x-x_0))$ vanishes for $x=x_0$ and $x=x_0+\pi/\sqrt{\delta}$, and  $\cos(\sqrt{\delta}(x-x_0))$ evaluated in these points becomes either $1$ or $-1$, which significantly simplifies corresponding definite integrals occurring in the sections above. 

\bigskip
Let us now provide a more detailed derivation of the expression $\|f^n_{\alpha,\beta}-\phi_n\|^2$ for odd $n \geq 3$ and $\alpha \geq n^2 \geq \beta$, see \eqref{eq:dist-fn-odd-a>b}.
Recall that each normalized Fu\v{c}\'ik eigenfunction $f^n_{\alpha,\beta}$  with odd $n$ has $\frac{n+1}{2}$ positive bumps and $\frac{n-1}{2}$ negative bumps. 
Therefore, using the formulas from above, we deduce that
\begin{align}
	\notag
	\|f^n_{\alpha,\beta}-\phi_n\|^2
	&=
	\sum_{i=0}^{\frac{n+1}{2}-1}\int_{il}^{il+l_1}
	\left(
	\frac{\sqrt{\beta}}{\sqrt{\alpha}}\sin(\sqrt{\alpha}(x-il))
	-
	\sin(nx)
	\right)^2 \mathrm{d}x \\
	\notag
	&+ \sum_{i=0}^{\frac{n-1}{2}-1}
	\int_{il+l_1}^{(i+1)l}
	\left(
	\sin(\sqrt{\beta}(x-il-l_1))
	+
	\sin(nx)
	\right)^2 \mathrm{d}x \\
	\notag
	&= 
	\frac{\pi}2 + \frac{n+1}{4}\frac{\beta}{\alpha}l_1
	+
	\frac{n-1}{4} l_2 \\
	\notag
	&+ 4\sqrt{\beta}\frac{\alpha-\beta}{(n^2-\alpha)(n^2-\beta)}\sum_{i=0}^{\frac{n-1}{2}-1}\left[\sin(n(il+l_1))+\sin(n(i+1)l)\right]\\
	\notag
	&= 
	\frac{\pi}2 + \frac{n+1}{4}\frac{\beta}{\alpha}l_1
	+
	\frac{n-1}{4} l_2 \\
	\notag
	&+ 8\sqrt{\beta}\frac{\alpha-\beta}{(n^2-\alpha)(n^2-\beta)}
	\frac{\sin\left(\frac{(n-1)nl}{8}\right)\cos\left(\frac{nl_2}2\right)}
	{\sin\left(\frac{nl}2\right)}
	\sin\left(\frac{(n-1)nl}{8}
	+
	\frac{nl_1}{2}\right)
	\\
	\notag
	&= \pi - \pi\frac{(n+1)(\sqrt{\alpha}-1)}{\sqrt{\alpha}(2\sqrt{\alpha}-(n+1))^2}(\sqrt{\alpha}-n) \\
	\notag
	&+ \frac{16(n-1)\sqrt{\alpha}^3}{(2\sqrt{\alpha}-(n+1))}\frac{(\sqrt{\alpha}-1)}{(n+\sqrt{\alpha})(n+1)(n^2-3n\sqrt{\alpha}+n+\sqrt{\alpha})}\times\\
	\label{eq:estim1}
	&\times \frac{\cos\left(\frac{\pi}2\frac{n}{\sqrt{\alpha}}\right)\cos\left(\frac{\pi}2\frac{n(2\sqrt{\alpha}-n-1)}{(n-1)\sqrt{\alpha}}\right)}{(\sqrt{\alpha}-n)\sin\left(\pi\frac{n(\sqrt{\alpha}-1)}{(n-1)\sqrt{\alpha}}\right)}.
\end{align}
Here, we used the summation formula
\begin{equation}\label{eq:sum}
\sum_{i=0}^{k-1}\sin(ci+d)=\frac{\sin\left(k\frac{c}2\right)\sin\left((k-1)\frac{c}2+d\right)}{\sin\left(\frac{c}2\right)}.
\end{equation}

Notice that the arguments of the last cosine and sine in \eqref{eq:estim1} satisfy
$$
\frac{\pi}{2} \leq \frac{\pi}2\frac{n(2\sqrt{\alpha}-n-1)}{(n-1)\sqrt{\alpha}} \leq \frac{3\pi}{2}
\quad \text{and} \quad
\pi \leq \pi\frac{n(\sqrt{\alpha}-1)}{(n-1)\sqrt{\alpha}} \leq \frac{3\pi}{2}.
$$
That is, the cosine and sine of these arguments are negative.
To make it easier to control the total sign in \eqref{eq:estim1}, we apply the formulas 
$$
\cos(x) = - \cos(\pi-x)
\quad\text{and}\quad
\sin(x) = - \sin(x-\pi).
$$ 
This gives the expression \eqref{eq:dist-fn-odd-a>b}.

\medspace
The derivation of the expression \eqref{eq:speon} for the scalar product $\langle f^n_{\alpha,\beta}, \phi_m \rangle$ for even $n$ and $m$ with $n>m$ follows similar steps as above:
\begin{align*}
\langle f^n_{\alpha,\beta},\phi_m\rangle &= \sum_{k=0}^{\frac{n}{2}-1}\frac{\sqrt{\beta}}{\sqrt{\alpha}}\int_{kl}^{kl+l_1}\sin(\sqrt{\alpha}(x-kl))\sin(mx)\,\mathrm{d}x \\
&- \sum_{k=0}^{\frac{n}{2}-1}\int_{kl+l_1}^{(k+1)l}\sin(\sqrt{\beta}(x-kl-l_1))\sin(mx)\,\mathrm{d}x \\
&=\sqrt{\beta}\left(\frac{1}{\alpha-m^2}-\frac{1}{\beta-m^2}\right)\sum_{k=0}^{\frac{n}{2}-1}\sin(mkl)+\sin(m(kl+l_1)),
\end{align*}
and we conclude by applying the summation formula \eqref{eq:sum}.

Let us now provide details on the derivation of the upper bound for $\|f^n_{\alpha,\beta}-\phi_n\|^2$ with odd $n \geq 3$ and $\alpha \geq n^2 \geq \beta$, given by \eqref{eq:asymp-fn-odd-1}.
Reducing the terms on the right-hand side of the first inequality in \eqref{eq:asymp-fn-odd-1} to a common denominator, we arrive at
\begin{align*}	
	\|f^n_{\alpha,\beta}-\phi_n\|^2
	&\leq
	\frac{\pi (\sqrt{\alpha}-n)^2}{144 (n-1) \alpha (2\sqrt{\alpha}-n-1)^2(n+\sqrt{\alpha})((3n-1)\sqrt{\alpha}-n(n+1))}
	\\
	&\times
	\bigg(
	144 \sqrt{\alpha} (n^2-1) (2 \sqrt{\alpha} + 4 \sqrt{\alpha^3} + n + n^2 - 2 \sqrt{\alpha} n (n+2) + \alpha (5n-7)) 
	\\
	&+ 
	48\pi^2 \alpha (\sqrt{\alpha}-1) (2 \sqrt{\alpha} - n - 1) (n^2+1)
	\\
	&-\pi^4 (\sqrt{\alpha}-1) (\sqrt{\alpha} - n)^2 (2 \sqrt{\alpha} - n - 1) (n+1)^2
	\bigg).
\end{align*}
Using now the following simple estimates:
\begin{align*}
&(2\sqrt{\alpha}-n-1)^2 
\geq 
\frac{(n-1)^2}{n^2}\alpha,
\quad
(n+\sqrt{\alpha})((3n-1)\sqrt{\alpha}-n(n+1))
\geq 
4n^2(n-1),\\
&n^2-1 \leq n^2+1, 
\quad 
n + n^2 - 2 \sqrt{\alpha} n (n+2) + \alpha (5n-7)
\leq
\alpha (5n-7) \leq 5\sqrt{\alpha^3},\\
&(\sqrt{\alpha}-1) (2 \sqrt{\alpha} - n - 1) \leq 2\alpha,
\quad
-\pi^4 (\sqrt{\alpha}-1) (\sqrt{\alpha} - n)^2 (2 \sqrt{\alpha} - n - 1) (n+1)^2 \leq 0,
\end{align*}
we get 
$$
\|f^n_{\alpha,\beta}-\phi_n\|^2
\leq
\frac{\pi (\sqrt{\alpha}-n)^2 (n^2+1) (288 + (1296 + 96 \pi^2) \alpha)}
{576 (n-1)^4 \alpha}.
$$
Finally, recalling that $\sqrt{\alpha} \geq n \geq 3$ and roughly estimating
$$
\frac{288 + (1296 + 96 \pi^2) \alpha}{576}
\leq 
4\alpha,
$$
we obtain
$$
\|f^n_{\alpha,\beta}-\phi_n\|^2
\leq
\frac{4\pi(n^2+1)(\sqrt{\alpha}-n)^2}{(n-1)^4} 
\leq
\frac{4\pi n^2 (n^2+1)}{(n-1)^4} 
\left(\frac{\sqrt{\alpha}}{n}-1\right)^2
\leq
23\pi\left(\frac{\sqrt{\alpha}}{n}-1\right)^2.
$$

For the case $\beta > n^2 > \alpha$, we use a similar procedure as above to estimate the upper bound for  $\|f^n_{\alpha,\beta}-\phi_n\|^2$ with odd $n \geq 3$ given by \eqref{eq:asymp:dist-fn-odd-b>a}.
We convert the right-hand side of the first inequality in \eqref{eq:asymp:dist-fn-odd-b>a} to a common denominator, to get
\begin{align*}
\|f^n_{\alpha,\beta}-\phi_n\|^2
&\leq
\frac{\pi (\sqrt{\beta}-n)^2}{144 (n+1)\beta {(2\sqrt{\beta}-n+1)^2} {(n+\sqrt{\beta})((3n+1)\sqrt{\beta}-n(n-1))}}
\\
&\times
\bigg(
144 \sqrt{\beta} {(n^2-1)} (2 \sqrt{\beta} + 4 \sqrt{\beta^3} + n - n^2 - 2 \sqrt{\beta}n(n-2) + \beta (5n+7)) 
\\
&+ 
48 \pi^2 \beta (\sqrt{\beta}+1) (2\sqrt{\beta} - n +1) (n^2+1)  
\\
&-  
{\pi^4(\sqrt{\beta}+1) (\sqrt{\beta} - 
	n)^2 (2\sqrt{\beta} -n + 1) (n-1)^2}
\bigg).
\end{align*}
With the simple estimates
\begin{align*}
	&(2\sqrt{\beta}-n+1)^2 
	\geq 
	\frac{(n+1)^2}{n^2}\beta,
	\quad
	(n+\sqrt{\beta})((3n+1)\sqrt{\beta}-n(n-1)) 
	\geq 4 n^2 (n+1),\\
	&n^2-1 \leq n^2 +1,
	\quad 
	n - n^2 - 2 \sqrt{\beta}n(n-2) + \beta (5n+7) 
	\leq \beta (5n+7) \leq 8 \sqrt{\beta^3},\\
	&(\sqrt{\beta}+1) (2\sqrt{\beta} - n +1) \leq 2\beta,
	\quad
	-\pi^4(\sqrt{\beta}+1) (\sqrt{\beta} - 
	n)^2 (2\sqrt{\beta} -n + 1) (n-1)^2 \leq 0,
\end{align*}
we obtain
\begin{align*}
\|f^n_{\alpha,\beta}-\phi_n\|^2
&\leq
\frac{\pi (\sqrt{\beta}-n)^2 (n^2+1)
{(288 + (1728 + 96 \pi^2) \beta)}}
{576 (n+1)^4\beta}.
\end{align*}
Finally, estimating
$$
\frac{288 + (1728 + 96 \pi^2) \beta}{576}
\leq 
5\beta,
$$
we obtain
$$
\|f^n_{\alpha,\beta}-\phi_n\|^2
\leq
\frac{5\pi(n^2+1)(\sqrt{\beta}-n)^2}{(n+1)^4} 
\leq
\frac{5\pi n^2 (n^2+1)}{(n+1)^4} 
\left(\frac{\sqrt{\beta}}{n}-1\right)^2
\leq
5\pi\left(\frac{\sqrt{\beta}}{n}-1\right)^2.
$$

\section{}\label{sec:appendix2}

In this section, we calculate the norms of the operators $T_k$ defined by \eqref{eq:TK}. 
First, let us show that each $T_k$ is linear. 
Indeed, taking any $g,h \in L^2(0,\pi)$, we have
$$
T_k(g+h)(x) 
= 
(g+h)^*\left(\frac{kx}{2}\right)
=
(-1)^\kappa (g+h)\left(\frac{kx}{2}-\pi \kappa\right)
$$
for $\pi \kappa \leq \frac{kx}{2} \leq \pi (\kappa+1)$, $\kappa \in \mathbb{N}_0$.
At the same time, for such $x$ we get
\begin{align*}
	(-1)^\kappa (g+h)\left(\frac{kx}{2}-\pi \kappa\right)
	&=
	(-1)^\kappa g\left(\frac{kx}{2}-\pi \kappa\right)
	+
	(-1)^\kappa h\left(\frac{kx}{2}-\pi \kappa\right)
	\\
	&=
	g^*\left(\frac{kx}{2}\right)
	+
	h^*\left(\frac{kx}{2}\right)
	=
	T_k g(x) + T_k h(x),
\end{align*}
which concludes the linearity.

Assume now that $k$ is even, i.e., $k=2m$, $m \geq 1$.
For any $g \in L^2(0,\pi)$ we obtain
\begin{align*}
	\int_0^\pi \left(T_{2m} g(x)\right)^2  \mathrm{d}x 
	&=
	\int_0^\pi \left(g^*(mx)\right)^2  \mathrm{d}x 
	=
	\frac{1}{m} \int_0^{\pi m} \left(g^*(x)\right)^2  \mathrm{d}x 
	\\
	&=
	\frac{1}{m} \sum_{\kappa=0}^{m-1} \int_{\pi \kappa}^{\pi (\kappa+1)} g^2(x-\pi \kappa) \, \mathrm{d}x 
	=
	\frac{1}{m} \sum_{\kappa=0}^{m-1} \int_0^\pi g^2(x)  \,\mathrm{d}x  
	=
	\int_0^\pi g^2(x)  \,\mathrm{d}x .
\end{align*}
Therefore, we get
$$
\|T_{2m}\|_* 
= 
\sup_{g \in L^2(0,\pi)\setminus\{0\}} \frac{\|T_{2m} g\|}{\|g\|} 
=1.
$$
Assume that $k$ is odd, i.e., $k=2m+1$, $m \geq 0$.
We have
\begin{align*}
	&\int_0^\pi \left(T_{2m+1} g(x)\right)^2  \mathrm{d}x 
	=
	\int_0^\pi \left(g^*\left(\frac{(2m+1)x}{2}\right)\right)^2  \mathrm{d}x 
	=
	\frac{2}{2m+1} \int_0^{\pi \left(m+\frac{1}{2}\right)} \left(g^*(x)\right)^2  \mathrm{d}x 
	\\
	&=
	\frac{2}{2m+1} \sum_{\kappa=0}^{m-1} \int_{\pi \kappa}^{\pi (\kappa+1)} g^2(x-\pi \kappa)  \,\mathrm{d}x 
	+
	\frac{2}{2m+1} \int_{\pi m}^{\pi \left(m+\frac{1}{2}\right)}  g^2(x-\pi m)  \,\mathrm{d}x 
	\\
	&=
	\frac{2m}{2m+1} \int_0^\pi g^2(x)  \,\mathrm{d}x  
	+
	\frac{2}{2m+1} \int_{0}^{\frac{\pi}{2}}  g^2(x)  \,\mathrm{d}x 
	\leq
	\frac{2m+2}{2m+1} \int_0^\pi g^2(x)  \,\mathrm{d}x .
\end{align*}
Notice that this estimate is sharp since equality holds for any $g$ with the support on  $(0,\pi/2)$.
Thus, we deduce that
$$
\|T_{2m+1}\|_* 
= 
\sup_{g \in L^2(0,\pi)\setminus\{0\}} \frac{\|T_{2m+1} g\|}{\|g\|} 
= 
\sqrt{\frac{2m+2}{2m+1}}.
$$


\begin{thebibliography}{99}

\bibitem{agran}
M. S. Agranovich, 
\textit{On series with respect to root vectors of operators associated with forms having symmetric principal part}, 
Functional Analysis and Its Applications, 28(3), 151-167, 1994.
\href{https://doi.org/10.1007/BF01078449}{\nolinkurl{DOI:10.1007/BF01078449}}

\bibitem{bary2}
N. Bary, \textit{Sur le syst\`{e}mes complets de fonctions othogonales}, Recueil math\'ematique de la Soci\'et\'e math\'ematique de Moscou [Matematicheskii Sbornik], 14(56), 1-2, 51-108, 1944.
\url{http://www.mathnet.ru/php/archive.phtml?wshow=paper&jrnid=sm&paperid=6184}

\bibitem{BB}
F. Baustian and V. Bobkov,
\textit{On asymptotic behavior of Dirichlet inverse},
International Journal of Number Theory, 16(6), 1337-1354, 2020.
\href{https://doi.org/10.1142/S1793042120500700}{\nolinkurl{DOI:10.1142/S1793042120500700}}

\bibitem{bind}
P. Binding, L. Boulton, J. \v{C}epi\v{c}ka, P. Dr\'{a}bek, and P. Girg, P., \textit{Basis properties of eigenfunctions of the $p$-Laplacian}, Proceedings of the American Mathematical Society, 134(12), 3487-3494, 2006.
\href{https://doi.org/10.1090/S0002-9939-06-08001-4}{\nolinkurl{DOI:10.1090/S0002-9939-06-08001-4}}

\bibitem{BL}
L. Boulton and G. Lord, 
\textit{Approximation properties of the $q$-sine bases}, 
Proceedings of the Royal Society A: Mathematical, Physical and Engineering Sciences, 467(2133), 2690-2711, 2011.
\href{https://doi.org/10.1098/rspa.2010.0486}{\nolinkurl{DOI:10.1098/rspa.2010.0486}}


\bibitem{BM}
L. Boulton and H. Melkonian, 
\textit{A multi-term basis criterion for families of dilated periodic functions}, 
Zeitschrift f\"ur Analysis und ihre Anwendungen, 38(1), 107-124, 2019.
\href{https://doi.org/10.4171/ZAA/1630}{\nolinkurl{DOI:10.4171/ZAA/1630}}

\bibitem{BE}
P.~J. Bushell and D.~E. Edmunds, 
\textit{Remarks on generalized trigonometric functions}, 
The Rocky Mountain Journal of Mathematics, 42(1), 25-57, 2012.
\url{https://www.jstor.org/stable/44240034}


\bibitem{cac}
N.~P. C\'ac, 
\textit{On nontrivial solutions of a Dirichlet problem whose jumping nonlinearity crosses a multiple eigenvalue}, 
Journal of Differential Equations, 80(2), 379-404, 1989.
\href{https://doi.org/10.1016/0022-0396(89)90090-9}{\nolinkurl{DOI:10.1016/0022-0396(89)90090-9}}

\bibitem{castro}
A. Castro and C. Chang, 
\textit{A variational characterization of the Fucik spectrum and applications}, 
Revista Colombiana de Matemáticas, 44(1), 23-40, 2010.
\url{https://revistas.unal.edu.co/index.php/recolma/article/view/28591}

\bibitem{cues}
M. Cuesta, D. de Figueiredo, and J. P. Gossez, \textit{The beginning of the Fu\v{c}ik spectrum for the $p$-Laplacian}, 
Journal of Differential Equations, 159(1), 212-238, 1999.
\href{https://doi.org/10.1006/jdeq.1999.3645U}{\nolinkurl{DOI:10.1006/jdeq.1999.3645}}


\bibitem{dancer}
E. N. Dancer, \textit{On the Dirichlet problem for weakly non-linear elliptic partial differential equations}, 
Proceedings of the Royal Society of Edinburgh Section A: Mathematics, 76(4), 283-300, 1977.
\href{https://doi.org/10.1017/S0308210500019648}{\nolinkurl{DOI:10.1017/S0308210500019648}}

\bibitem{dancer2}
E. N. Dancer,
\textit{Remarks on jumping nonlinearities}, 
in J. Escher and G. Simonett (eds) \textit{Topics in nonlinear analysis. Progress in Nonlinear Differential Equations and Their Applications}, pp. 101-116, 1999. Birkhäuser, Basel.
\href{https://doi.org/10.1007/978-3-0348-8765-6_7}{\nolinkurl{DOI:10.1007/978-3-0348-8765-6_7}}

\bibitem{djak}
P. Djakov and B. Mityagin, 
\textit{Criteria for existence of Riesz bases consisting of root functions of Hill and $1$D Dirac operators}, 
Journal of Functional Analysis, 263(8), 2300-2332, 2012.
\href{https://doi.org/10.1016/j.jfa.2012.07.003}{\nolinkurl{DOI:10.1016/j.jfa.2012.07.003}}

\bibitem{drabek}
P. Dr\'abek and S. B. Robinson, \textit{A new and extended variational characterization of the Fu\v{c}\'ik Spectrum with application to nonresonance and resonance problems}, 
Calculus of Variations and Partial Differential Equations, 57(1), 1, 2018.
\href{https://doi.org/10.1007/s00526-017-1276-8}{\nolinkurl{DOI:10.1007/s00526-017-1276-8}}

\bibitem{duff}
R.~J. Duffin and J.~J. Eachus, \textit{Some Notes on an Expansion Theorem of Paley and Wiener},
Bulletin of the American Mathematical Society, 48(12), 850-855, 1942.
\url{https://projecteuclid.org/euclid.bams/1183504861}


\bibitem{EGL}
D.~E. Edmunds, P. Gurka, and J. Lang, 
\textit{Properties of generalized trigonometric functions}, 
Journal of Approximation Theory, 164(1), 47-56, 2012.
\href{https://doi.org/10.1016/j.jat.2011.09.004}{\nolinkurl{DOI:10.1016/j.jat.2011.09.004}}


\bibitem{fucik}
S. Fu\v{c}\'ik, \textit{Boundary value problems with jumping nonlinearities}, 
\v{C}asopis pro p\v{e}stov\'an\'i matematiky, 101(1), 69-87, 1976.
\url{http://dml.cz/dmlcz/108683}


\bibitem{HLS}
H. Hedenmalm, P. Lindqvist, and K. Seip, \textit{A Hilbert space of Dirichlet series and systems of dilated functions in $L^2(0,1)$}, 
Duke Mathematical Journal, 86(1), 1-37, 1997.
\href{https://doi.org/10.1215/S0012-7094-97-08601-4}{\nolinkurl{DOI:10.1215/S0012-7094-97-08601-4}}


\bibitem{HLS2}
H. Hedenmalm, P. Lindqvist, and K. Seip, \textit{Addendum to ``A Hilbert space of Dirichlet series and systems of dilated functions in $L^2(0,1)$''}. Duke Mathematical Journal, 99(1), 175-178, 1999.
\href{https://doi.org/10.1215/S0012-7094-99-09907-6}{\nolinkurl{DOI:10.1215/S0012-7094-99-09907-6}}


\bibitem{holubova}
G. Holubov\'a and P. Ne\v{c}esal, 
\textit{The Fu\v{c}\'ik spectrum: exploring the bridge between discrete and continuous world}, 
in S. Pinelas, M. Chipot, and Z. Dosla (eds)
\textit{Differential and Difference Equations with Applications}, pp. 421-428. Springer,  2013.
\href{https://doi.org/10.1007/978-1-4614-7333-6_36}{\nolinkurl{DOI:10.1007/978-1-4614-7333-6_36}}

\bibitem{kato}
T. Kato, \textit{Perturbation Theory for Linear Operators}, Springer, 1980.
\href{https://doi.org/10.1007/978-3-642-66282-9}{\nolinkurl{DOI:10.1007/978-3-642-66282-9}}

\bibitem{kats}
V. E. Katsnel'son, 
\textit{Conditions under which systems of eigenvectors of some classes of operators form a basis}, 
Functional Analysis and Its Applications, 1(2), 122-132, 1967.
\href{https://doi.org/10.1007/BF01076084}{\nolinkurl{DOI:10.1007/BF01076084}}


\bibitem{li}
C. Li, S. Li, Z. Liu, and J. Pan, 
\textit{On the Fuc\'ik spectrum}, 
Journal of Differential Equations, 244(10), 2498-2528, 2008.
\href{https://doi.org/10.1016/j.jde.2008.02.021}{\nolinkurl{DOI:10.1016/j.jde.2008.02.021}}

\bibitem{locker}
J. Locker, 
\textit{Spectral theory of non-self-adjoint two-point differential operators}, 
American Mathematical Society, 2000.
\href{https://doi.org/10.1090/surv/073}{\nolinkurl{DOI:10.1090/surv/073}}

\bibitem{loos}
I. Looseov\'a and P. Ne\v{c}esal, 
\textit{The Fu\v{c}\'ik spectrum of the discrete Dirichlet operator}, 
Linear Algebra and its Applications, 553, 58-103, 2018.
\href{https://doi.org/10.1016/j.laa.2018.04.017}{\nolinkurl{DOI:10.1016/j.laa.2018.04.017}}

\bibitem{markus}
A. S. Markus, 
\textit{Introduction to the spectral theory of polynomial operator pencils}. 
American Mathematical Society, 1988.

\bibitem{paley}
R.~C. Paley and N. Wiener, \textit{Fourier Transforms in the Complex Domain}, American Mathematical Society, 1934.
\url{https://bookstore.ams.org/coll-19}

\bibitem{perera}
K. Perera and M. Schechter, \textit{The Fucik spectrum and critical groups}, 
Proceedings of the American Mathematical Society, 129(8), 2301-2308, 2001.
\href{https://doi.org/10.1090/S0002-9939-01-05968-8}{\nolinkurl{DOI:10.1090/S0002-9939-01-05968-8}}

\bibitem{pinsal}
J. P. Pinasco and A. M. Salort, \textit{Asymptotic behavior of the curves in the Fu\v{c}\'ik spectrum}, 
Communications in Contemporary Mathematics, 19(04), 1650039, 2017.
\href{https://doi.org/10.1142/S0219199716500395}{\nolinkurl{DOI:10.1142/S0219199716500395}}

\bibitem{RS}
M. Reed and B. Simon, 
\textit{Methods of modern mathematical physics I: Functional analysis}, 
Academic Press, 1981.
\href{https://doi.org/10.1016/B978-0-12-585001-8.X5001-6}{\nolinkurl{DOI:10.1016/B978-0-12-585001-8.X5001-6}}

\bibitem{rynne}
B. P. Rynne, 
\textit{The Fucik spectrum of general Sturm-Liouville problems}, 
Journal of Differential Equations, 161(1), 87-109, 2000.
\href{https://doi.org/10.1006/jdeq.1999.3661}{\nolinkurl{DOI:10.1006/jdeq.1999.3661}}

\bibitem{skal1}
A. A. Shkalikov, 
\textit{Perturbations of self-adjoint and normal operators with discrete spectrum}, 
Russian Mathematical Surveys, 71(5), 907, 2016.
\href{https://doi.org/10.1070/RM9740}{\nolinkurl{DOI:10.1070/RM9740}}

\bibitem{schech2}
M. Schechter, 
\textit{The Fucik spectrum}, 
Indiana University Mathematics Journal, 43(4), 1139-1157, 1994.
\href{https://doi.org/10.1512/iumj.1994.43.43050}{\nolinkurl{DOI:10.1512/iumj.1994.43.43050}}

\bibitem{sing}
I. Singer, \textit{Bases in Banach Spaces I}, Springer, 1970.
\url{https://www.springer.com/gp/book/9783642516351}

\bibitem{wintner2}
A. Wintner, \textit{On T\"opler's Wave Analysis}, 
American Journal of Mathematics, 69(4), 758-768, 1947.
\href{https://doi.org/10.2307/2371797}{\nolinkurl{DOI:10.2307/2371797}}

\bibitem{young}
R.~M. Young, \textit{An Introduction to Nonharmonic Fourier Series}, Academic Press, 1980. 


\end{thebibliography}
\end{document}